\documentclass[11pt]{article}
\usepackage{amsfonts,mathrsfs,amssymb,amsthm,mathptm}
\usepackage{amsmath,amscd}
\usepackage{mathptm,pslatex}
\usepackage{color}
\usepackage[dvips]{graphicx}
\usepackage[all]{xy}
\usepackage{graphicx}
\usepackage{fancyhdr}
\usepackage{url}

\begin{document}
\newtheorem{defi}{Definition}[section]
\newtheorem{exam}[defi]{Example}
\newtheorem{prop}[defi]{Proposition}
\newtheorem{theorem}[defi]{Theorem}
\newtheorem{lem}[defi]{Lemma}
\newtheorem{coro}[defi]{Corollary}
\theoremstyle{definition}
\newtheorem{rem}[defi]{Remark}
\newtheorem{ques}[defi]{Question}

\newcommand{\add}{{\rm add}}
\newcommand{\con}{{\rm con}}
\newcommand{\gd}{{\rm gl.dim}}
\newcommand{\sd}{{\rm st.dim}}
\newcommand{\sr}{{\rm sr}}
\newcommand{\dm}{{\rm dom.dim}}
\newcommand{\cdm}{{\rm codomdim}}
\newcommand{\tdim}{{\rm dim}}
\newcommand{\E}{{\rm E}}
\newcommand{\Mor}{{\rm Morph}}
\newcommand{\End}{{\rm End}}
\newcommand{\rank}{{\rm rank}}
\newcommand{\PSL}{{\rm PSL}}
\newcommand{\GL}{{\rm GL}}
\newcommand{\ind}{{\rm ind}}
\newcommand{\rsd}{{\rm res.dim}}
\newcommand{\rd} {{\rm rd}}
\newcommand{\ol}{\overline}
\newcommand{\overpr}{$\hfill\square$}
\newcommand{\rad}{{\rm rad}}
\newcommand{\soc}{{\rm soc}}
\renewcommand{\top}{{\rm top}}
\newcommand{\pd}{{\rm pdim}}
\newcommand{\id}{{\rm idim}}
\newcommand{\fld}{{\rm fdim}}
\newcommand{\Fac}{{\rm Fac}}
\newcommand{\Gen}{{\rm Gen}}
\newcommand{\fd} {{\rm fin.dim}}
\newcommand{\Fd} {{\rm Fin.dim}}
\newcommand{\Pf}[1]{{\mathscr P}^{<\infty}(#1)}
\newcommand{\DTr}{{\rm DTr}}
\newcommand{\cpx}[1]{#1^{\bullet}}
\newcommand{\D}[1]{{\mathscr D}(#1)}
\newcommand{\Dz}[1]{{\mathscr D}^+(#1)}
\newcommand{\Df}[1]{{\mathscr D}^-(#1)}
\newcommand{\Db}[1]{{\mathscr D}^b(#1)}
\newcommand{\C}[1]{{\mathscr C}(#1)}
\newcommand{\Cz}[1]{{\mathscr C}^+(#1)}
\newcommand{\Cf}[1]{{\mathscr C}^-(#1)}
\newcommand{\Cb}[1]{{\mathscr C}^b(#1)}
\newcommand{\Dc}[1]{{\mathscr D}^c(#1)}
\newcommand{\K}[1]{{\mathscr K}(#1)}
\newcommand{\Kz}[1]{{\mathscr K}^+(#1)}
\newcommand{\Kf}[1]{{\mathscr  K}^-(#1)}
\newcommand{\Kb}[1]{{\mathscr K}^b(#1)}
\newcommand{\DF}[1]{{\mathscr D}_F(#1)}

\newcommand{\Kac}[1]{{\mathscr K}_{\rm ac}(#1)}
\newcommand{\Keac}[1]{{\mathscr K}_{\mbox{\rm e-ac}}(#1)}

\newcommand{\modcat}{\ensuremath{\mbox{{\rm -mod}}}}
\newcommand{\cmodcat}{\ensuremath{\mbox{{\rm -comod}}}}
\newcommand{\Modcat}{\ensuremath{\mbox{{\rm -Mod}}}}
\newcommand{\ires}{\ensuremath{\mbox{{\rm ires}}}}
\newcommand{\Stires}{\ensuremath{\mbox{{\rm Stires}}}}
\newcommand{\Stpres}{\ensuremath{\mbox{{\rm Stpres}}}}
\newcommand{\Spec}{{\rm Spec}}

\newcommand{\stmc}[1]{#1\mbox{{\rm -{\underline{mod}}}}}
\newcommand{\Stmc}[1]{#1\mbox{{\rm -{\underline{Mod}}}}}
\newcommand{\prj}[1]{#1\mbox{{\rm -proj}}}
\newcommand{\inj}[1]{#1\mbox{{\rm -inj}}}
\newcommand{\Prj}[1]{#1\mbox{{\rm -Proj}}}
\newcommand{\Inj}[1]{#1\mbox{{\rm -Inj}}}
\newcommand{\PI}[1]{#1\mbox{{\rm -Prinj}}}
\newcommand{\GP}[1]{#1\mbox{{\rm -GProj}}}
\newcommand{\GI}[1]{#1\mbox{{\rm -GInj}}}
\newcommand{\gp}[1]{#1\mbox{{\rm -Gproj}}}
\newcommand{\gi}[1]{#1\mbox{{\rm -Ginj}}}

\newcommand{\opp}{^{\rm op}}
\newcommand{\otimesL}{\otimes^{\rm\mathbb L}}
\newcommand{\rHom}{{\rm\mathbb R}{\rm Hom}\,}
\newcommand{\pdim}{\pd}
\newcommand{\Hom}{{\rm Hom}}
\newcommand{\Coker}{{\rm Coker}}
\newcommand{ \Ker  }{{\rm Ker}}
\newcommand{ \Cone }{{\rm Con}}
\newcommand{ \Img  }{{\rm Im}}
\newcommand{\Ext}{{\rm Ext}}
\newcommand{\StHom}{{\rm \underline{Hom}}}
\newcommand{\StEnd}{{\rm \underline{End}}}
\newcommand{\KK}{I\!\!K}
\newcommand{\gm}{{\rm _{\Gamma_M}}}
\newcommand{\gmr}{{\rm _{\Gamma_M^R}}}

\def\demo{{\bf Proof}\hskip10pt}

\def\g{\gamma} \def\d{\delta} \def\a{\alpha}
\def\s{\sigma}  \def\om{\omega}  \def\ld{\lambda}
\def\D{\Delta}
\def\si{\Sigma} \def\O{\Omega}
\def\G{\Gamma} \def\GG{{\cal G}} \def\XX{{\cal X}} \def\MM{{\cal M}} \def\lrr{\lg r\rg } \def\ogg{\overline {\GG}}
\def\og{\overline G} \def\oh{\overline H} \def\oc{\overline C}
 \def\oq{\overline Q}

 \def\oa{\overline A}
  \def\ob{\overline B}  \def\ol{\overline L} \def\om{\overline M}
\def\on{\overline N} \def\op{\overline P} \def\os{\overline S}
\def\ot{\overline T} \def\ok{\overline K} \def\ov{\overline V}
\def\od{\overline D} \def\oi{\overline I} \def\oj{\overline J}
\def\o1{\overline 1} \def\olh{\overline h}

\def\o{\overline}   \def\olr{\overline r} \def\oll{\overline \ell} \def\olt{\overline t }

\def\di{\bigm|} \def\lg{\langle} \def\rg{\rangle}

\def\vez{\varepsilon}\def\bz{\bigoplus}  \def\sz {\oplus}
\def\epa{\xrightarrow} \def\inja{\hookrightarrow}

\newcommand{\lra}{\longrightarrow}
\newcommand{\llra}{\longleftarrow}
\newcommand{\lraf}[1]{\stackrel{#1}{\lra}}
\newcommand{\llaf}[1]{\stackrel{#1}{\llra}}
\newcommand{\ra}{\rightarrow}
\newcommand{\dk}{{\rm dim_{_{k}}}}

\newcommand{\holim}{{\rm Holim}}
\newcommand{\hocolim}{{\rm Hocolim}}
\newcommand{\colim}{{\rm colim\, }}
\newcommand{\limt}{{\rm lim\, }}
\newcommand{\Add}{{\rm Add }}
\newcommand{\Prod}{{\rm Prod }}
\newcommand{\pres}{\ensuremath{\mbox{{\rm pres}}}}
\newcommand{\app}{{\rm app }}
\newcommand{\Tor}{{\rm Tor}}
\newcommand{\Cogen}{{\rm Cogen}}
\newcommand{\Tria}{{\rm Tria}}
\newcommand{\Loc}{{\rm Loc}}
\newcommand{\Coloc}{{\rm Coloc}}
\newcommand{\tria}{{\rm tria}}
\newcommand{\Con}{{\rm Con}}
\newcommand{\Thick}{{\rm Thick}}
\newcommand{\thick}{{\rm thick}}
\newcommand{\Sum}{{\rm Sum}}
\def\Mon{\hbox{\rm Mon}}
\def\Aut{\hbox{\rm Aut}}
\newcommand{\PGL}{{\rm PGL}}
\def\Syl{\hbox{\rm Syl}}
\def\char{\hbox{\rm \,char\,}}
{\Large \bf
\begin{center}
A classification of regular maps with Euler characteristic $-p^4$ for a prime $p\geq 5$
\end{center}}

\medskip
\centerline{\textbf{Xiaogang Li$^*$} and \textbf{Yao Tian} }

\medskip

\renewcommand{\thefootnote}{\alph{footnote}}
\setcounter{footnote}{-1} \footnote{ $^*$ Corresponding author.
Email: 2200501002@cnu.edu.cn.}
\renewcommand{\thefootnote}{\alph{footnote}}
\setcounter{footnote}{-1}
\renewcommand{\thefootnote}{\alph{footnote}}
\setcounter{footnote}{-1}

\begin{abstract}
A map is a cellular decomposition of a closed surface. In the framework of classifying all regular maps by their supporting surface, it is an open problem to find all closed surfaces that support no regular maps. Classification of regular maps on surfaces with Euler characteristic $-p, -p^2, -p^3, -2p,$ and $-3p$ has already been done by several authors in a series of papers, which also show that surfaces with these Euler characteristic support no regular maps if the corresponding prime $p$ satisfies certain conditions. In this paper, assuming that $p\geq 5$ is a prime and $i\geq 4$, we show that the order of a Sylow $p$-subgroup of a regular map with Euler characteristic $-p^i$ is bounded by $p^{i-1}$ unless $p\in \{5, 7, 13\}$, and we show the existence of a normal $p$-subgroup for these regular maps whenever a Sylow $p$-subgroup has order at least $\sqrt{p^i}$, laying a solid foundation for using an inductive method to completely characterize regular maps of Euler characteristic $-p^i$. Based on this, we classify all regular maps with Euler characteristic $-p^4$ for a prime $p\geq 5$ in terms of reduced presentations of their automorphism groups. Consequently, a closed surface with Euler characteristic $-p^4$ supports no regular maps if and only if $p\notin \{2,3,5,7,13\}$.
\end{abstract}

\section{Introduction}
A  map is a cellular decomposition of a closed surface. In other words, a {\it map} $\MM$ is an embedding (without crossings) of graphs (or multigraphs) $\cal{G}$ into a closed surface $\cal{S}$, such that connected components of $\cal{S}\setminus \cal{G}$ are simply connected. Given a map $\MM$ on a closed surface $\cal{S}$ with underlying graph $\cal{G}$, the connected components of $\cal{S}\setminus \cal{G}$ are called {\it faces} of the map $\MM$, and the vertices and edges of $\cal{G}$ are called respective {\it vertices} and {\it edges} of $\MM$. 
By an {\it automorphism} of a map $\MM$ we mean an  automorphism of the {\it underlying graph} ${\cal G}$
which extends to a  self-homeomorphism of the surface.
These automorphisms  form a subgroup $\Aut({\MM})$ of the automorphism  group $\Aut({\cal G})$ of $\cal G$.
It is well-known that  $\Aut({\MM})$ acts semi-regularly on the set of all flags (in most cases, which are incident vertex-edge-face triples).
If the action is regular, then we call the map as well as the corresponding embedding  {\it regular}.  The map is called orientable (non-orientable) if the supporting surface is orientable (non-orientable). The Euler characteristic of a regular map is referred to that of its supporting surface.
\vskip 2mm
Regular maps are closely related to other branches of mathematics such as group theory, hyperbolic geometry and complex curves and so on. During the past forty years, there have been fruitful research on regular maps, see  \cite{BS,CCDKNW,CDNS, GNSS,JJ,JS} for an overview.
\vskip 2mm

There are three approaches to the study of regular maps: classification by surface (the topological approach); classification by underlying graph (the graph-theoretical approach); classification by automorphism group (the group-theoretical approach). In this paper, we shall focus on the topological approach-classification by surface. For a very long time, classification of regular maps has been known on only a finite list of surfaces \cite{MC, MCD}. The first breakthrough was the classification of regular maps with Euler characteristic $-p$ by A. D'azevedo, R. Nedela and J. \v{S}ir\'{a}\v{n} in \cite{ARJ}, which starts the classification of regular maps with Euler characteristic $-p^n$ for $p$ a prime and $n$ a positive integer. Several years later, M. Conder, P. Poto\v{c}nik, and J. \v{S}ir\'{a}\v{n} gave a classification of regular maps with almost Sylow-cyclic groups in \cite{MPJ} and successfully applied this to classify regular maps with Euler characteristic $-p^2$. Their work emphasizes the notion of almost Sylow-cyclic groups and sheds light on research in this direction. After that, regular maps with Euler characteristic $-2p, -3p,-p^3$ and orientably-regular maps with Euler characteristic $-2p^2$ have also been classified in two decades, see \cite{CNJ,CST,Ma,TL} for an overview. The techniques for all these classifications come mainly from group theory. By \cite{ARJ}, one knows that a closed surface with Euler characteristic $-p$ supports no regular maps if $p\geq 17$ and $p\equiv 1({\rm mod}~ 12)$, and it is known from \cite{MPJ} that a closed surface with Euler characteristic $-p^2$ supports no regular maps provided $p\notin \{3,7\}$, and \cite{TL} proves that a closed surface with Euler characteristic $-p^3$ supports no regular maps if $p\equiv 1({\rm mod}~ 12)$. Motivated by these works, we raise the following two questions:
\begin{ques}\label{q}
What is the set $\mathcal{A}$ of all positive integers $i$ such that there exist only finitely many primes $p$ for which closed surfaces with Euler characteristic $-p^i$ support regular maps?
\end{ques}
\begin{ques}\label{q2}
What is the set $\mathcal{B}$ of all positive integers $i$ such that all but finitely many primes $p$ satisfy that closed surfaces with Euler characteristic $-p^i$ support only solvable regular maps (regular maps with solvable automorphism groups)?
\end{ques}
For Questions \ref{q} and \ref{q2}, we note that $\mathcal{A}\cap \mathcal{B}=\emptyset$, $2\in \mathcal{A}$ is implied by \cite[Theorem 6.1]{MPJ} and $1,3\in \mathcal{B}$ are implied by the main theorems in \cite{ARJ} and \cite{TL}, respectively. Understanding the sets $\mathcal{A}$ and $\mathcal{B}$ is important for the following two reasons:
\vskip 2mm
$(1)$ It enhances the classification of regular maps with Euler characteristic $-p^i$ for a prime $p$;
\vskip 2mm
$(2)$ It enlarges the list of surfaces that do not support regular maps. 
\vskip 2mm
In this paper, we classify regular maps with Euler characteristic $-p^4$ for a prime $p\geq 5$. 
\begin{theorem}\label{main}
A regular map $\MM$ is of Euler characteristic $-p^4$ for a prime $p\geq 5$ if and only if one of the following holds
\begin{enumerate}
\item[\rm (1)] $p=5$, $\MM$ is isomorphic to $\MM_1(1,3,4)$ or its dual;
\item[\rm (2)] $p=7$, $\MM$ is isomorphic to $\MM_4(1,1,4)$ or $\MM_6(1,1,3)$ or their dual;
\item[\rm (3)] $p=13$, $\MM$ is isomorphic to $\MM_7(1,7,6)$ or its dual.
\end{enumerate}
where $\MM_1(1,3,4), \MM_4(1,1,4), \MM_6(1,1,3)$ and $\MM_7(1,7,6)$ are defined in Section 4.
\end{theorem}
Note that these exist regular maps of Euler characteristic $-2^4$ (all of them are orientable) and regular maps of Euler characteristic $-3^4$ (all of them are non-orientable) by the list in \cite{MC}. As a consequence of Theorem \ref{main} and this fact, we have the following result, which proves $4\in \mathcal{A}$. 
\begin{coro}\label{main1}
Let $p$ be a prime. Then a closed surface with Euler characteristic $-p^4$ supports no regular maps if and only if $p\in \{2, 3, 5, 7, 13\}$.
\end{coro}

\begin{rem}
During our writing this paper, we noticed that \cite{MNJ} provides a rough characterization of regular maps $\MM=\MM(G;r,t,\ell)$ with Euler characteristic $-p^k$ (where $p$ is an odd prime and $k$ a positive integer) in terms of the types of $\MM$ and the quotient groups $G/N$ over $N$ with $N$ being either the largest normal subgroup of odd order or a normal $p$-subgroup, and Corollary \ref{main1} is implied by \cite[Corollary 1.1]{MNJ}. But our works differs from \cite{MNJ} and has independent significance due to the following three reasons:

$(1)$ Methodology-the methodology used in this paper inherits from our previous work $\cite{TL}$ which is different from that of \cite{MNJ}, where \cite[Lemma 3.5]{TL} is one of the cornerstones of our method;

$(2)$ Aim-compared to the rough characterizations given in \cite{MNJ}, our goal (realized in \cite{TL} and this paper) is to characterize regular maps with Euler characteristic $-p^i$ for $p$ a prime and $i$ a fixed number using an inductive method:specifically, characterizing regular maps of Euler characteristic $-p^i$ when all regular maps of Euler characteristic $-p^j$ with $j<i$ are already known. We believe that regular maps with Euler characteristic $-p^i$ for $p$ a prime and $i$ a relatively small integer can be completely classified in this way;

$(3)$ Result-\cite[Corollary 1.1]{MNJ} provides a list of types and rough characterization of the automorphism groups for regular maps with Euler characteristic $-p^4$ for $p\in \{5,7,13\}$, we give explicitly the reduced representations for the automorphism groups of these regular maps and we reveal the importance of studying irreducible representations of certain 2-dimensional projective special (general) linear groups in characterizing regular maps with Euler characteristic $-p^i$.
\end{rem}
After this introductory section,  some notations, terminologies, known  results and  basic theory of regular maps  will be introduced in  Section 2, and Theorem~\ref{main} will be proved in Section 3. In Section 4, we will characterize regular maps with Euler characteristic $-5^4, -7^4$ and $-13^4$ and confirm the existence of them.

\section{Preliminaries}
In this section, we shall give  some  notations and terminologies,  a brief introduction to regular maps and some known results used in this paper.

\subsection{Notations and terminologies }

For a finite set $\Omega$, the cardinality of $\Omega$ is denoted by $|\Omega|$.
For a given prime $p$ and an arbitrary integer $m$, the $p$-part and $p'$-part of $m$ refer to the largest $p$-power dividing $m$ and the largest divisor of $m$ coprime to $p$, respectively.

\vskip 2mm
Let $G$ be a finite group. For $H<G$, $H\unlhd G$ and $H~{\rm char} G$, we mean that $H$ is a proper subgroup of $G$, a normal subgroup of $G$ and a characteristic subgroup of $G$, respectively. For $g\in G$ and $H\leq G$, we use $|g|$ and $|G:H|$ to denote the order of $g$ and the index of $H$ in $G$, respectively. For a prime $p$, we denote the set of Sylow $p$-subgroups of $G$ (if $p\nmid |G|$, then a Sylow $p$-subgroup of $G$ is identified as the trivial subgroup $\{1\}$) and the number of Sylow $p$-subgroups of $G$ by $\Syl_p(G)$ and $n_p(G)$, respectively. The largest normal $p$-subgroup of $G$ is denoted by $O_{p}(G)$ and the largest normal subgroup of $G$ of odd order is denoted by $O_{2'}(G)$. The centre of $G$, the Fitting subgroup of $G$ (the product of all nilpotent normal  subgroups of $G$) is denoted by $F(G)$ and the commutator subgroup of $G$ are denoted by $Z(G), F(G)$ and $G'$, respectively.
\vskip 2mm

For finite groups $G$ and $A$, if there is a homomorphism $\phi: A\rightarrow \Aut(G)$, then the semiproduct of $G$ and $A$ with respect to $\phi$ is denote by $G\rtimes_\phi A$.
\vskip 2mm
By $\GL(n,q)$, $\PGL(n,q)$ and $\PSL(n,q)$, we denote the general linear groups, projective general linear groups and projective special linear groups of dimension $n$ over the finite field $F_q$, respectively. By $\mathbb{Z}_m$ and $\mathbb{D}_{2n}$ ($n\geq 2$), we denote the cyclic group of order $m$ and dihedral group of order $2n$, respectively.
\vskip 2mm
We introduce an important  terminology which was defined in  \cite{MPJ} as follows.  A group $G$ is called {\it almost Sylow-cyclic} if  all of its Sylow subgroups of odd order  are cyclic and  all of its  Sylow-$2$ subgroup  are either trivial or contain a cyclic subgroup of index $2$.

\subsection{Regular maps}

In this subsection, we give a description of maps and regular maps in the combinatorial way. We begin with the following definition.
\begin{defi}\label{map1}
{\rm For a given finite set $F$ and three fixed-point-free involutory permutations $r, t, \ell $ on $F$, a quadruple
$\MM=\MM(F; r, t, \ell )$ is called a {\it combinatorial map} if they satisfy the following two conditions:

(1)\ $t$ commutes with $\ell$; 

(2)\ the group $\lg r,t, \ell \rg $ acts transitively on $F$.}
\end{defi}

For a combinatorial map $\MM=\MM(F; r,t, \ell )$,  $F$  is called the {\it flag} set,
$r, t, \ell$ are called {\it rotary, transversal and longitudinal involution,} respectively.
The group $\lg r, t, \ell \rg $ is called the {\it monodromy group} of $\MM$, denoted by $\Mon(\MM)$.
We define the {\it vertices, edges} and {\it face-boundaries} of $\MM$ to be the
orbits of the subgroups $\lg r,t\rg$, $\lg t, \ell \rg $ and $\lg r, \ell \rg $ respectively.
The incidence in $\MM $ is defined by a nontrivial set intersection.

Given two maps $\MM _1=\MM(F_1; r_1, t_1, \ell _1)$ and $\MM_2=\MM_2(F_2; r_2, t_2, \ell _2)$, a bijection $\phi$ from $F_1$ to
$F_2$ is called a {\it map isomorphism} if $\phi r_1=r_2\phi$, $\phi t_1=t_2\phi$ and $\phi \ell _1=\ell _2\phi$. In particular,
if $\MM _1=\MM _2=\MM$, then $\phi$ is called an {\it automorphism} of $\MM$.
The set $\Aut(\MM)$ of all automorphisms of $\MM$ forms a group  which is called the {\it automorphism group} of $\MM$. By the definition of map isomorphism, we have $\Aut(\MM )=C_{S_F}(\Mon (\MM))$,  the centralizer of $\Mon(\MM)$ in
$S_F$. It follows from the transitivity of $\Mon(\MM )$ on $F$ that $\Aut(\MM)$ acts semi-regularly on $F$.  If the action is regular, then we call $\MM$
{\it regular map}.

Let $\MM=\MM(F; r,t, \ell )$ be a regular map and fix a flag $\alpha\in F$. Then for every $\beta\in F$, there is a unique $g_{\beta}\in \Aut(\MM)$ such that $\beta=\alpha^{g_{\beta}}$ and the mapping
$$\varphi: F\rightarrow \Aut(\MM),~\beta\mapsto g_{\beta}$$
is a bijection. Let $R(\Aut(\MM))$ be the right regular representation of $\Aut(\MM)$. Then the mapping
 $$\phi: \Aut(\MM)\rightarrow R(\Aut(\MM)),~h\mapsto R(h), \forall~ h\in \Aut(\MM)$$ is a group isomorphism.
Since $$\varphi(\beta^h)=\varphi(\alpha^{g_{\beta}h})
=g_{\beta}h=g_{\beta}^{R(h)}=\varphi(\beta)^{\phi(h)}$$
for every $\beta\in F$ and $h\in \Aut(\MM)$, it follows that the two permutation groups $\Aut(\MM)\le{\rm Sym}(F)$ and $R(\Aut(\MM))\leq {\rm Sym}(\Aut(\MM))$ are permutation isomorphic.
In this sense, we identify $F$ with the group $G=\langle r, t, \ell\rangle$ and the rotary, transversal and longitudinal involution on $F$ with the left regular representation of $r$, $t$ and $\ell$, respectively. Thus the monodromy
(automorphism) group of $\MM$ is the left (right) regular representation of $G$. Usually, we identify $R(G)$ with $G$ and call $\MM:=\MM(G; r, t, \ell )$ an {\it algebraic map}. It is obvious that $\MM^*:=\MM(G; r, \ell, t)$ is also an algebraic map. We call it the \emph{dual} of $\MM$.
It is straightforward to check that the even-word subgroup $\lg tr, r\ell \rg $ of $G$ has index at most 2.
If the index is 2, then one may fix an orientation for $\MM $ and so $\MM $ is said to be {\it orientable}.  If the index is 1, then $\MM$ is said to be {\it non-orientable}. The {\it type} of $\MM$ is defined as $\{|rt|,|r\ell|\}$, and the \emph{Euler characteristic} $\chi_{\MM}$ of $\MM$ is meant the Euler characteristic of the supporting surface of $\MM$, that is, $\chi_{\MM}:=V+F-E$. In particular, if we set $x:=|rt|$ and $y:=|r\ell|$, then we obtain the following formula for $\chi_{\MM}$:
$$
(*)\quad\quad\chi_{\MM}=-\frac{|G|(xy-2x-2y)}{4xy}.$$ Clearly, $\chi_{\MM}=\chi_{\MM^*}$. Whenever $xy-2x-2y$ is positive (equivalently, $\chi_{\MM}$ is negative), the well-known Hurwitz bound says that $\frac{xy}{xy-2x-2y}\leq 21$ and the equality is attained for $\{x,y\}=\{3,7\}$. In particular, $|G|\leq -84\chi_{\MM}$ in this case. The well-known Euler-Poincar\'{e} formula relates the genera $g_{\cal{S}}$ of a surface $\cal{S}$ to its Euler characteristic $\chi_{\cal{S}}$ as follows:
$$\chi_{\cal{S}}=\begin{cases}2-2g_{\cal{S}} & \mbox{ if }~\cal{S}~\mbox{is orientable},\\ 2-g_{\cal{S}} & \mbox{ if }~\cal{S}~\mbox{is non-orientable}.\end{cases}$$Thus only non-orientable surfaces can have odd Euler characteristic.

Following \cite[Section 1]{TL}, for two coprime odd integers $j$ and $k$, we denote by $\MM(j,k)$ the regular map $\MM(G;r,t,\ell)$ over the group $G=\mathbb{D}_{2j}\times \mathbb{D}_{2k}$ of type $\{2j,2k\}$; for any integer $m$ with $m\equiv 3({\rm mod}~6)$, we denote by $\MM(m)$ the regular map $\MM(G;r,t,\ell)$ over the group $G=\mathbb{D}_4\rtimes_\phi \mathbb{D}_{2m}$ (where $\phi$ is a surjective homomorphism from $\mathbb{D}_{2m}$ onto $\Aut(\mathbb{D}_4)$) of  type $\{4,m\}$. 

\subsection{Some Known Results}

\begin{prop}[{\cite[Chap.1, Theorem 6.11]{SUZ}}]
\label{nc}

Let $H$ be a subgroup of a group $G$.
Then $C_G(H)$ is a normal subgroup of
$N_G(H)$ and the quotient $N_G(H)/C_G(H)$ is isomorphic
to a subgroup of $\Aut (H)$.
\end{prop}

\begin{prop}[{\cite[Chap.1, Theorem 1.16]{ISA}}]
\label{num}
Suppose that $G$ is a finite group such that $n_p(G)>1$, and choose distinct Sylow $p$-subgroups $S$ and $T$ such that the order $|S\cap T|$ is as large as possible. Then $n_p(G)\equiv 1\pmod{|S:S\cap T|}$.
\end{prop}

\begin{prop}[{\cite[Chap.5, Theorem 5.3]{ISA}}]
\label{is}
Let $G$ be a finite group and $p$ be a prime divisor of $|Z(G)\cap G'|$.
Then a Sylow-$p$ subgroup of $G$ is non-abelian.
\end{prop}

\begin{prop}[{\cite[Lemma 3.2]{MPJ}}]
\label{alm}
Let $\MM(G;r,t,\ell)$ be a regular map. If $p$ is a prime divisor of $|G|$ coprime
to the Euler characteristic of $\MM$, then a Sylow $p$-subgroup of $G$ is cyclic (if $p$ is odd) or dihedral (if $p=2$).
\end{prop}

\begin{prop}[{\cite[Theorem 5.4]{MPJ}}]
\label{non-alm}

Let $\MM(G;r,t,\ell)$ be a regular map. If $G$ is an insolvable almost Sylow-cyclic group, then there exist a prime $q$ and an integer $d$ coprime to $q(q^2-1)$ such that $G\cong \PSL(2,q)$ or $\mathbb{Z}_d\rtimes_\phi \PGL(2,q)$, where $\phi$ is the homomorphism from $\PGL(2,q)$ onto $\langle\sigma_{\mathbb{Z}_d}\rangle$ where $\sigma_{\mathbb{Z}_d}$ is automorphism of $\mathbb{Z}_d$ sending $i$ to $-i$ for every $i\in\mathbb{Z}_d$.
\end{prop}

\begin{prop}[{\cite[Theorem 1]{GW}}]
\label{dihedral}
Let $G$ be a finite group with dihedral Sylow $2$-subgroups. Then $G/O_{2'}(G)$ is isomorphic to either an almost simple group with socle $\PSL(2, q)$ where $q$ is odd,
or $A_7$, or a Sylow $2$-subgroup of $G$.
 \end{prop}

\begin{prop}[{\cite[Table 1]{MPJ}}]
\label{solvable}
Let $\MM(G;r,t,\ell)$ be a non-orientable regular map with $G$ a solvable almost Sylow-cyclic group. Then we have one of the following:
\begin{enumerate}
\item[{\rm (i)}] $\MM$ has type $\{2,n\}$ and $G\cong \mathbb{D}_{2n}$, where $n\geq 4$ is even;
\item[{\rm (ii)}] $\MM$ or $\MM^*\cong \MM(m,n)$, where $2\nmid mn$, $1\neq m<n$ and $\gcd(m,n)=1$;
\item[{\rm (iii)}] $\MM$ or $\MM^*\cong \MM(m)$, where $m\equiv 3({\rm mod}~6)$,
where $\MM(m)$ and  $\MM(m,n)$ are defined in Section 2.2.
\end{enumerate}
 \end{prop}
 
\begin{prop}[{\cite[Theorem II.7.3a, p.187]{Hup}}]
\label{HU}
Let $G$ be an irreducible subgroup of a Singer cycle S of $\GL(n,p)$. Then the centralizer in $\GL(n,p)$ of $G$ is $S$, and the normalizer in $\GL(n,p)$ of $G$ is the semiproduct of $S$ and a cyclic group $C$ of order $n$.
 \end{prop} 

\begin{prop}[{\cite[Lemma 3.5]{TL}}]
\label{normal}
Let $\MM(G;r,t,\ell)$ be a regular map of Euler characteristic $-p^k$ for some odd prime $p$ and positive integer $k\geq 3$ and let $P\in \Syl_p(G)$. If $O_p(G)=1$, then $O_{2'}(G)$ is a cyclic group of order not divided by $p$ and one of the followings holds
\begin{enumerate}
\item[\rm (1)] $P=1$, $\MM(G;r,t,\ell)\cong \MM(m,n)$ or $\MM(m,n)^*$, where $2\nmid mn$, $\gcd(m,n)=1$, $1\neq m<n$
and $mn-m-n=p^k$;
\item[\rm (2)] $P=1$, $\MM(G;r,t,\ell)\cong \MM(m)$ or $\MM(m)^*$ , where $m\equiv 3\pmod 6$ and $m-4=p^k$;
\item[\rm (3)] $P$ is cyclic and $G\cong \mathbb{Z}_d\rtimes_\phi \PGL(2,q)$, where $q$ is an odd prime, $\gcd(d, q(q^2-1))=1$ and $\phi$ is the homomorphism from $\PGL(2,q)$ onto $\langle \sigma_{\mathbb{Z}_d}\rangle$;
\item[\rm (4)] $P$ is elementary abelian of rank at least $2$, $G/O_{2'}(G)$ is isomorphic to an almost simple group $N$ with socle $\PSL(2,r)$ and
$|N:\PSL(2,r)|\leq 2$, where $r=|P|$.
\end{enumerate}
In particular, if $P$ is neither cyclic nor elementary abelian, then $O_p(G)>1$.
\end{prop}

\begin{prop}[{\cite[Lemma 3.6]{TL}}]
\label{cop}
Let $G$ be either $\PSL(2,q)$ or $\mathbb{Z}_d\rtimes_\phi\PGL(2,q)$ for some odd prime $q\geq 5$, where $\phi$ is the homomorphism from $\PGL(2,q)$ onto $\langle \sigma_{\mathbb{Z}_d}\rangle$. If $G=\langle a,b\rangle$
for two elements $a$ and $b$ such that any odd prime divisor of $|G|$ divides $|a||b|$, then $\gcd(|a|,|b|)=1$.
\end{prop}
We point out that \cite[Theorem 2.2(3)]{ARJ} excludes the case-a regular map of type $\{4, 6\}$ with Euler characteristic $-5$ and automorphism group $\PGL(2,5)$. The following is a consequence of this fact and \cite[Theorem 2.2]{ARJ}.
\begin{prop}
\label{cla}
Let $\MM(G;r,t,\ell)$ be a non-orientable regular map with Euler characteristic $-p$ for some prime $p\geq 5$. Suppose that either $G$ is insolvable or $G$ does not have a normal Sylow $p$-subgroup. Then one of the following holds
 \begin{enumerate}
\item[{\rm (i)}] $p=5$, $\MM$ has type $\{4,6\}$ and $G\cong \PGL(2,5)$;
\item[{\rm (ii)}] $p=7$, $\MM$ has type $\{3,8\}$ and $G\cong \PGL(2,7)$;
\item[{\rm (iii)}] $p=13$, $\MM$ has type $\{3,7\}$ and $G\cong \PSL(2,13)$.
\end{enumerate}
 \end{prop}
The following result is an immediate corollary of \cite[Theorem 6.1]{MPJ}.
\begin{prop}
\label{Eule}
Let $\MM$ be a non-orientable regular map with Euler characteristic $-p^2$ for some prime $p\geq 5$. Then $\MM$ is of type $\{3,13\}$, $\Aut(\MM)\cong \PSL(2,13)$ and $p=7$.
\end{prop} 
The following result is an immediate consequence of \cite[Theorem 1.1]{TL}
\begin{prop}\label{p^3}
Let $\MM(G;r,t,\ell)$ be a regular map with Euler characteristic $-p^3$ for some odd prime $p\geq 5$. Then $G$ has trivial Sylow $p$-subgroup and is solvable. 
\end{prop}

\section{Proof of Theorem \ref{main}}
The proof of Theorem \ref{main} is purely group-theoretical. Before we proceed, we fix a condition, mention a few basic facts, and prove two general results on regular maps with Euler characteristic $-p^i$ for a prime $p\geq 5$ and $i\geq 4$. 
\vskip 2mm
For convenience, we fix the following condition:
$$(\star)~~\MM(G; r, t, \ell)~~\mbox{is a regular map with Euler characteristic}~-p^{i}~
\mbox{where}~p\geq 5~\mbox{is a prime and }~ i\geq 4.$$
Assume condition $(\star)$. Then we have the following three basic facts:

$(1)$ $G=\lg rt, r\ell\rg$;
\vskip 2mm
$(2)$ $\lg rt\rg \cap \lg r\ell\rg =1$;
\vskip 2mm
$(3)$ $|rt|\geq 3$ and $|r\ell|\geq 3$.
\vskip 2mm
\noindent These facts will be used throughout this section without any further reference.
\vskip 2mm
The following general lemma will be crucial in our later analysis. 
\begin{theorem}\label{divide}
Assume condition $(\star)$. Then $p^i\nmid |G|$ unless $p\in \{5,7,13\}$.
\end{theorem}
\demo We abbreviate $\MM(G;r,t,\ell)$ by $\MM$ and let $x$ and $y$ be the order of $rt$ and $r\ell$, respectively. Suppose for the contrary that $p^i\mid |G|$. By the formula $(*)$, we deduce that $k(x,y)=\frac{xy}{xy-2x-2y}$ is a positive integer. It is known from the proof of Theorem 3.4 in \cite{MPJ} that there are finitely many pairs $\{x,y\}$ such that $k(x,y)$ is a positive integer, which are listed as follows:
\begin{table}[h]
\centering
\caption{List of \{x,y\} and corresponding integer value of k(x,y)}
$
\begin{tabular}{c|c|c|c|c|c}
                 \hline
                 \{x,y\} & k(x,y) & \{x,y\} & k(x,y) & \{x,y\} & k(x,y) \\
                 \hline
                 \{3,7\} & 21 & \{3,24\} & 4 & \{5,5\} & 5\\

                 \{3,8\} & 12 & \{4,5\} & 10 & \{5,20\} & 2\\

                 \{3,9\} & 9 &  \{4,6\} & 6 & \{6,6\} & 3\\

                 \{3,12\}& 6 & \{4,8\}& 4 & \{6,12\}& 2 \\

                 \{3,15\}& 5 & \{4,12\}& 3 & \{8,8\}& 2\\
                 \hline
                
               \end{tabular}$
               \end{table}

We claim that $p\nmid xy$. Suppose for the contrary that $p\mid xy$. Let $P$ be a Sylow $p$-subgroup of $G$. Then, checking the list, we find that $p=5$ and $\{x,y\}=\{3,15\},\{4,5\},\{5,5\},\{5,20\}$ or $p=7$ and $\{x,y\}=\{3,7\}$. Assume $p=5$. By computing the order of $G$ with formula $(*)$, we can exclude the case $\{x,y\}=\{3,15\}$ as $3\nmid |G|$ in this case. For the remaining cases, by Sylow's theorems we obtain $P\unlhd G$. Consider the quotient group $\og:=G/P$. We find that $\og\neq \langle \olr\olt, \olr\oll\rangle$. Thus $G\neq \langle rt,r\ell\rangle$, a contradiction. Hence $p=7$ and $\{x,y\}=\{3,7\}$. This together with Sylow's theorems implies that $G$ has a normal Sylow $7$-subgroup, say $R$. Since $R$ and the quotient group $G/R$ are solvable, it follows that $G$ is solvable, too. In particular, $G'<G$. On the other hand, it follows from $G=\lg r, t, \ell\rg$  that the quotient group $G/G'$ is an abelian $2$-group. This together with $\{x,y\}=\{3,7\}$ implies that $rt, r\ell\in G'$. Thus $G=\langle rt,r\ell\rangle\leq G'<G$, a contradiction. The claim follows. We now proceed in the following two steps.
\vskip 2mm
{\it Step 1:  Show that $|P:O_p(G)|\leq p$.}
\vskip 2mm
If $P$ is normal, then $P=O_p(G)$ and so $|P:O_p(G)|=1$, as desired. Thus we may assume that $P$ is non-normal. In particular, $n_p(G)>1$, where $n_p(G)$ is the number of Sylow $p$-subgroups of $G$. Checking the list, we find that $1$ is the only divisor of $|G|$ that is congruent to $1$ modulo $p^2$. Thus $n_p(G)\not\equiv 1({\rm mod}~p^2)$. Applying Proposition \ref{num}, we know that there exists two distinct Sylow $p$-subgroups of $G$, say $R$ and $Q$, such that $|R:R\cap Q|=p$. We show that $R\cap Q\unlhd G$. To see it, with an observation on the list in the table, we find that $|G|$ has at most one divisor different from $1$ that is congruent to $1$ modulo $p$, except for the cases $p=5$ and $\{x,y\}=\{3,7\}, \{3, 8\}$ or $\{3, 9\}$. 

Suppose first that either $p\neq 5$, or $\{x,y\}\neq \{3,7\}, \{3, 8\}$ or $\{3, 9\}$. Set $N=N_G(R\cap Q)$. Then both $R$ and $Q$ are contained in $N$. In particular, $n_p(N)>1$. Note that $n_p(N)$ and $n_p(G)$ are divisors of $|G|$ that is congruent to $1$ modulo $p$. Thus the set of Sylow $p$-subgroups of $N$ is the same as that of $G$. Let $M$ be the normal subgroup of $G$ which is generated by all of its Sylow $p$-subgroups. Then $\lg R, Q\rg \leq M\leq N$. Since $N=N_G(R\cap Q)$, it follows that $R\cap Q\unlhd M$. In particular, $R\cap Q=O_p(M)$. As $O_p(M)\char M\unlhd G$, we have $R\cap Q=O_p(M)\unlhd G$. Thus it follows from $O_p(M)\leq O_p(G)$ that $|P:O_p(G)|\leq |P:O_p(M)|=p$, as desired.

Next, we consider the cases $p=5$ and $\{x,y\}=\{3,7\}, \{3, 8\}$ or $\{3, 9\}$. Suppose $p=5$ and $\{x,y\}=\{3,7\}$. Then $|G|=2^2\cdot 3\cdot 5^4\cdot 7$. Since $G$ has dihedral Sylow $2$-subgroup by Proposition \ref{alm}, it follows from Proposition \ref{dihedral} that $G/O_{2'}(G)$ is isomorphic to $A_7$, a $2$-group or an almost simple group with socle $\PSL(2,q)$ for some odd $q$. The case $G/O_{2'}(G)\cong A_7$ can be excluded as the former one can be generated by its involutions but the latter one does not. Assume that $G$ is insolvable. Then it follows that $G/O_{2'}(G)$ is isomorphic to an almost simple group with socle $\PSL(2,q)$. In particular, $\frac{q(q^2-1)}{2}$ divides $|G|$. By comparing $\frac{q(q^2-1)}{2}$ with $|G|$, we find that $q=5$. Thus $|O_{2'}(G)|=5^3\cdot 7$. By Sylow's theorems, we deduce that $O_{2'}(G)$ has a normal Sylow $5$-subgroup, say $K$. Then $K\char O_{2'}(G)\unlhd G$ and therefore $K\unlhd G$. Thus $|P:O_5(G)|\leq |P:K|=5$, as desired. Assume that $G$ is solvable. Then it is easy to see that $G/O_{2'}(G)$ is a $2$-group ( the case $G/O_{2'}(G)\cong A_4$ can be excluded by the same reason as for $G/O_{2'}(G)\cong A_7$). In particular, $|O_{2'}(G)|=3\cdot 5^4\cdot 7$. Then a similar argument as in previous paragraph shows that $|O_5(O_{2'}(G))|\geq 5^3$. This implies that $|O_5(G)|\geq 5^3$. Thus $|P:O_5(G)|\leq 5$, as desired. Similarly, for the cases $p=5$ and $\{x,y\}=\{3, 8\}$ or $\{3, 9\}$, we have $|P:O_5(G)|\leq 5$, as desired.
\vskip 2mm
{\it Step 2:  Show the quotient map induced by $O_p(G)$ does not exist whenever $p\notin \{5, 7,13\}$.}
\vskip 2mm
Let $\overline{\MM}$ be the quotient map of $\MM$ induced by $O_p(G)$. If $P=O_p(G)$, then $\chi_{\overline{\MM}}=-1$, it is known from \cite{MCD} that such map does not exist.  So we may assume that $|P:O_p(G)|=p$. Then a Sylow $p$-subgroup of $G/O_p(G)$ is non-normal and $\overline{\MM}$ is a non-orientable map of Euler characteristic $-p$. According to Proposition \ref{cla}, there are three possibilities: $p=5, \{x, y\}=\{4,6\}$ and $G/O_p(G)\cong \PGL(2,5)$; $p=7, \{x,y\}=\{3,8\}$ and $G/O_p(G)\cong \PGL(2,7)$; $p=13$, $\{x,y\}=\{3,7\}$ and $G/O_p(G)\cong {\rm PSL}(2,13)$. The statement follows. \qed

\begin{rem}
The proof of \cite[Theorem 3.1]{TL} has a small gap in Step 1, as the cases $p=5$ and $\{x,y\}=\{3,7\}, \{3, 8\}$ or $\{3, 9\}$ must be considered independently. However, Step 1 of the revised proof corrects this gap.
\end{rem}

One of the obstacles in characterizing a regular map $\MM(G;r,t,\ell)$ with Euler characteristic $-p^i$ is that the structure of Sylow $p$-subgroups of $G$ is unknown. One way to overcome this problem is to passing from $G$ to its quotient group over a normal $p$-subgroup. The following lemma is crucial for our reductive method.
\begin{prop}\label{normal-sub}
Assume condition $(\star)$ and let $P\in \Syl_p(G)$. Then $O_p(G)>1$ whenever $|P|^2\geq p^i$.
\end{prop}
\demo Let $p^j$ be the order of $P$. We use contradiction, suppose that $O_p(G)=1$ and $j\geq \lceil\frac{i}{2}\rceil$, where $\lceil k\rceil$ for $k\in \mathbb{Q}$ means the smallest integer that is no less than $k$. Then by Proposition \ref{normal}, we have one of the following two cases:

\begin{enumerate}
\item[\rm (1)] $P$ is cyclic and $G\cong \mathbb{Z}_d\rtimes_\phi \PGL(2,q)$, where $q$ is an odd prime, $p\nmid d, \gcd(d, q(q^2-1))=1$ and $\phi$ is the homomorphism from $\PGL(2,q)$ onto $\langle \sigma_{\mathbb{Z}_d}\rangle$;
\item[\rm (2)] $P$ is elementary abelian of rank at least $2$, $G/O_{2'}(G)$ is isomorphic to an almost simple group $N$ with socle $\PSL(2,p^j)$ and
$|N:\PSL(2,p^j)|\leq 2$.
\end{enumerate}

Assume case (1). Then $P$ is cyclic. The assumption $p^{2j}\geq p^i\geq p^3$ implies that $j\geq 2$ and $p\neq q$. Therefore $p^j\mid q+1$ or $p^j\mid q-1$. Since $q+1$ and $q-1$ are even and $p^j$ is odd, it follows that $q+1\geq 2p^j$ in either cases. Thus $|G|=dq(q^2-1)\geq 2p^j(2p^j-1)(2p^j-2)$. But on the other, from the Hurwitz bound we have $|G|\leq 84p^i$. Hence $2p^j(2p^j-1)(2p^j-2)\leq 84p^i$. However, this together with $p\geq 5$ yields that $2\cdot 5^2(5^2-12)+1\leq 2p^j(p^j-12)+1\leq 0$, a contradiction.

Assume case (2). Then $$|G|=|O_{2'}(G)||G/O_{2'}(G)|\geq |G/O_{2'}(G)|\geq |\PSL(2,p^j)|= \frac{p^j(p^{2j}-1)}{2}.$$ This together with  $|G|\leq 84p^i$ implies that $\frac{p^j(p^{2j}-1)}{2}\leq 84p^i$. Let $x$ and $y$ denote $|rt|$ and $|r\ell|$, respectively, and let $x'$ and $y'$ denote the order of $|O_{2'}(G)rt|$ and $|O_{2'}(G)r\ell|$ in $G/O_{2'}(G)$, respectively. Then $x'\mid x$ and $y'\mid y$. By formula $(*)$, we see that, for any odd prime divisor $s$ of $|G|$ different from $p$, $|G|_s$ divides at least one of $x$ and $y$. Thus any odd prime divisors of $|G/O_{2'}(G)|$ different from $p$ divides one of $x'$ and $y'$. On the other hand, it follows from $G/O_{2'}(G)\cong N$ with $|N:\PSL(2,p^j)|\leq 2$ that $x'_{2'}$ (resp. $y'_{2'}$) divides one of $p, \frac{p^j-1}{2}$ and $\frac{p^j+1}{2}$. In particular, $x'_{2'}$ (resp. $y'_{2'}$) divides one of $p, (\frac{p^j-1}{2})_{2'}$ and $(\frac{p^j+1}{2})_{2'}$. This implies that $\{x'_{2'}, y'_{2'}\}=\{(\frac{p^j-1}{2})_{2'}, (\frac{p^j+1}{2})_{2'}\}$. Note that the assumption $p\geq 5$ is equivalent to either $p=5$ or $p\geq 7$. In what follows, we consider these two cases separately.

Suppose $p\geq 7$. We first show that $j=2$ and $i=4$. Suppose for the contrary that either $j\geq 3$ or $i\geq 5$. Then $j\geq 3$ in either cases, and from the inequality $\frac{p^j(p^{2j}-1)}{2}\leq 84p^i$ we obtain $$\frac{7^3(p^i-1)}{2}\leq \frac{p^3(p^i-1)}{2}\leq \frac{p^j(p^{2j}-1)}{2}\leq 84p^i.$$ Hence $\frac{7^3(p^i-1)}{2}\leq 84p^{i}$. This implies $\frac{1}{p^{i}}\geq \frac{175}{343}$, a contradiction to the fact that $\frac{1}{p^{i}}\leq \frac{1}{7^5}<\frac{175}{343}$. Now, substituting $j=i=2$ into the inequality $\frac{p^j(p^{2j}-1)}{2}\leq 84p^{2i}$, we obtain $p^2(p^4-1)\leq 168p^4$. Solving it, we have $p\in \{7, 11, 13\}$ under our assumption. For the three cases $p=7, p=11$ and $p=13$, we deduce that $$\frac{1}{x'_{2'}}+ \frac{1}{y'_{2'}}\leq \max\{\frac{1}{3}+\frac{1}{25}, \frac{1}{15}+\frac{1}{61}, \frac{1}{17}+\frac{1}{21}\}=\frac{28}{75}.$$Hence $$\frac{1}{x}+\frac{1}{y}\leq \frac{1}{x'}+\frac{1}{y'}\leq \frac{1}{x'_{2'}}+\frac{1}{y'_{2'}}\leq\frac{28}{75}.$$But then we obtain $$24p^4=\frac{49(p^4-\frac{p^4}{49})}{2}\leq\frac{49(p^4-1)}{2}\leq |\PSL(2,p^2)|\leq |G|=\frac{4p^4}{1-2(\frac{1}{x}+\frac{1}{y})}\leq \frac{4p^4}{1-\frac{56}{75}}<24p^4.$$ This is clearly a contradiction. 

Now we are left with the case $p=5$. From the inequality $\frac{p^j(p^{2j}-1)}{2}\leq 84p^i$ and the assumption $j\geq \lceil \frac{i}{2}\rceil\geq 2$, it is easy to see that $i=4$ and $j\in \{2,3\}$. With a similar argument on $x$ and $y$ as for the case $p\geq 7$, we may also derive a contradiction. \qed
\vskip 2mm

The following result is crucial for our proof of Theorem \ref{main}.
\begin{prop}\label{sylow-p}
Let $\MM(G; r, t, \ell)$ be a regular map with Euler characteristic $-p^{4}$ for a prime
$p\geq 5$. Then $p^4\mid |G|$.
\end{prop}
\demo For convenience, we fix the following condition:
$$(\star\star)~~ \MM=\MM(G; r, t, \ell)~\mbox{is a regular map with}~ \chi_\MM=-p^{4}~\mbox{for some prime
}~p\geq 5~\mbox{and}~p^4\nmid |G|.$$ To prove the statement, it suffices to show that $\MM$ does not exist under the condition $(\star\star)$. Assume condition $(\star\star)$ and let $P\in \Syl_p(G)$. In the following two subsections, we shall derive two facts:
\begin{enumerate}
\item[{\rm (1)}]  $G$ is solvable and $P\unlhd G$ (Theorem~\ref{so-norm}); and
\item[{\rm (2)}]  $P$ is trivial (Theorem~\ref{cyclic}).
\end{enumerate}

By Proposition~\ref{alm}, we deduce that $G$ is a solvable almost Sylow-cyclic group (see Section 2 for its definition). Furthermore,  by Proposition \ref{solvable}, one of the following holds:
\begin{enumerate}
\item[\rm (1)] $\{x,y\}=\{2,n\}$, $n$ even, $G\cong \mathbb{D}_{2n}$;
\item[\rm (2)] $\MM\cong \MM(m)$ or $\MM(m)^*$ , where $m\equiv 3({\rm mod}~6)$,
\item[\rm (3)]  $\MM\cong \MM(m,n)$ or $\MM(m,n)^*$, where $2\nmid mn$, $1\neq m<n$ and $\gcd(m,n)=1$;

where  $\MM(m)$ and  $\MM(m,n)$ are defined in Section 2.
 \end{enumerate}

The first case  can be easily excluded, as $\chi_\MM=1\neq -p^4$. Assume case (2). Then $-p^4=\chi_\MM=4-m$. Thus $0\equiv m=p^4+4\equiv 2({\rm mod}~3)$, a contradiction. Assume case (3). Then $-p^4=\chi_\MM=m+n-mn$. Thus $2\equiv p^4+1=(m-1)(n-1)\equiv 0({\rm mod}~4)$, a contradiction.  \qed
\vskip 2mm
Now, we are in a position to prove Theorem \ref{main}.
\vskip 2mm
\noindent{\bf Proof of Theorem \ref{main}} 
Let $\MM(G; r, t, \ell)$ be a regular map with Euler characteristic $-p^{4}$ for some prime
$p\geq 5$. Then it follows from Theorem \ref{divide} and Proposition \ref{sylow-p} that $p\in \{5, 7, 13\}$. According to Theorem \ref{cha-exc}, the Theorem \ref{main} follows. \qed
\subsection{Solvability of $G$ and normality of $P$ }
In this subsection, we mainly prove:
\begin{theorem}\label{so-norm}
Assume condition $(\star\star)$ and let $P\in \Syl_p(G)$. Then $G$ is solvable and $P\unlhd G$.
\end{theorem}
A crucial step in the proof of Theorem \ref{so-norm} is to exclude the case that $G$ is an insolvable almost Sylow-cyclic group. The will be done in the following several lemmas.

\vskip 2mm
The following elementary lemma will be needed in later proof. It can be obtained by a case-by-case analysis, we omit its proof. 
\begin{lem}\label{int}
Let $j, k$ be two positive integers such that $j+k=4$, and let $p\geq 5$ be a prime. Then the set $\{\frac{2p^k+4p^j-4}{(2p^j-1)(2p^j-4)}, \frac{2p^k+4p^j+4}{2p^j(2p^j+1)},\frac{2p^k+4p^j-2}{(2p^j-2)(2p^j-3)}, \frac{2p^k+4p^j+2}{(2p^j-1)(2p^j+2)}\}$ contains no integers.
\end{lem}
Assume condition $(\star\star)$. The following lemma gives some characterization of  $G$ when $G$ is an insolvable almost Sylow-cyclic group.
\begin{lem}\label{notdi}
Assume condition $(\star\star)$ and suppose that $G$ is an insolvable almost Sylow-cyclic group. Then $G$ is isomorphic to either $\PSL(2,p)$ or $\mathbb{Z}_d\rtimes_\phi \PGL(2,p)$ and $p\nmid |rt||r\ell|$.
\end{lem}
\demo We abbreviate $\MM(G;r,t,\ell)$ by $\MM$ and let $\{x,y\}$ be the type of $\MM$. According to Proposition \ref{non-alm}, we know that $G$ is isomorphic to either $\PSL(2,q)$ or $\mathbb{Z}_d\rtimes_\phi {\rm PGL}(2,q)$ for some prime $q\geq 5$, and $d$ divides exactly one of $x$ and $y$ in the latter case.  Without loss of generality, we may assume that $d\mid x$. From formula $(*)$, we obtain the following equation:
$$
(\dag)~~~~~ (xy-2x-2y)\cdot q \cdot\frac{q+1}{2} \cdot\frac{q-1}{2}=\alpha \cdot\frac{x}{s}\cdot y\cdot p^4,
$$where $\alpha=2$ and $s=1$ if $G\cong \PSL(2,q)$, and $\alpha=1$ and $s=d$ if $G\cong\mathbb{Z}_d\rtimes_\phi {\rm PGL}(2,q)$.

We claim that $q\nmid xy$. Note that $q=p$ would follow immediately from Equation $(\dag)$ if the claim were proved. Suppose to the contrary that $q\mid xy$. Then $q\mid x$ or $q\mid y$. Assume that $q\mid \gcd(x,y)$. Then $q= \frac{x}{s}$, $q=y$. Substituting this into the Equation $(\dag)$, we have $(sq-2s-2) \frac{q+1}{2} \frac{q-1}{2}=\alpha p^4$. Since $\gcd(\frac{q+1}{2}, \frac{q-1}{2})=1$, we deduce from this equation that $2\in \{\frac{q+1}{2}, \frac{q-1}{2}\}$ and $p$ divides one of $\frac{q+1}{2}$ and $\frac{q-1}{2}$. Thus $q=5$ and $p=3$, a contradiction to our assumption $p\geq 5$. This implies that $q$ divides exactly one of $\frac{x}{s}$ and $y$. Hence, $\{q\}\subsetneq \{\frac{x}{s},y\}$. Let $h\in \{\frac{x}{s},y\}\backslash{q}$. By \cite[Table 2]{MPJ}, $h\mid \frac{q\pm 1}{2}$ if $G\cong \PSL(2,q)$, $h\mid {q\pm 1}$ but $h\nmid \frac{q\pm 1}{2}$ if $G\cong\mathbb{Z}_d\rtimes_\phi {\rm PGL}(2,q)$. In what follows, we consider the two cases $G\cong \PSL(2,q)$ and $G\cong\mathbb{Z}_d\rtimes_\phi {\rm PGL}(2,q)$, separately.
\vskip 2mm
{\it Case $(1)$: $G\cong \PSL(2,q)$.}
\vskip 2mm
In this case $s=1$, we may assume without loss of generality that $q\mid x$ but $q\nmid y$. Then $x=q$ and $y$ divides $\frac{q-1}{2}$ or $\frac{q+1}{2}$. In particular, $\gcd(x,y)=1$ and therefore $\gcd(x,xy-2x-2y)=\gcd(x,2)\leq 2$. We first consider the case $y\mid\frac{q-1}{2}$. Then from Equation $(\dag)$ we see that both $\frac{q+1}{2}$ and $xy-2x-2y$ divide $2p^4$. Note that $\frac{q+1}{2}\geq \frac{5+1}{2}=3$. It follows $p\mid \frac{q+1}{2}$, equivalently, $q\equiv -1({\rm mod}~p)$. Suppose $p\nmid xy-2x-2y$. Then $xy-2x-2y\leq 2$. Solving this inequality, we obtain $q=x=5$ and $y=4$. But then $y=4\nmid \frac{q-1}{2}=2$, a contradiction to our assumption. Thus $p\mid xy-2x-2y$. Substituting $x=q\equiv -1({\rm mod}~p)$ into this congruence equation, we obtain $xy-2x-2y\equiv -y +2 -2y({\rm mod}~p)$. On the other hand, it follows from the Equation $(\dag)$ that $y=\frac{q-1}{4}$ or $\frac{q-1}{2}$. Substituting $y$  in the congruence equation $-3y +2 \equiv 0({\rm mod}~p)$, we get $3q\equiv 11({\rm mod}~p)$ when $y=\frac{q-1}{4}$, and $3q\equiv 7({\rm mod}~p)$ when  $y=\frac{q-1}{2}$. This together with $q\equiv -1({\rm mod}~p)$ implies that $p=5$ when $y=\frac{q-1}{2}$, and $p=7$ when $y=\frac{q-1}{4}$. Substituting $p$ and $y$ into the Equation $(\dag)$, we obtain $2\cdot 7^4=\frac{q+1}{2}\cdot \frac{q^2-aq+2}{2}$, where $a=7$ when $p=5$, and $a=11$ when $p=7$. However, it is straightforward to see that no prime is the solution to this equation, a contradiction to the fact that $q$ is a prime.  Similarly, we obtain two congruence equations $q\equiv 1({\rm mod}~p)$ and $y +2 \equiv 0({\rm mod}~p)$ for the case $y\mid\frac{q+1}{2}$. Then $p=3$ when $y=\frac{q+1}{2}$, and $p=5$ when $y=\frac{q+1}{4}$. By our assumption, we have $p=5$ when $y=\frac{q+1}{4}$. Substituting $p$ and $y$ into the Equation $(\dag)$, we obtain $2\cdot 5^4=\frac{q-1}{2}\cdot \frac{q^2-9q-2}{2}$. Plainly, no prime is the solution to this equation, a contradiction to the fact that $q$ is a prime. Thus case (1) cannot occur.
\vskip 2mm
{\it Case $(2)$: $G\cong \mathbb{Z}_d\rtimes_\phi \PGL(2,q)$.}
\vskip 2mm
In this case, $\alpha=1$ and $s=d$. We have to consider four cases: $q\mid \frac{x}{s}$ and $y\mid q-1$;$q\mid \frac{x}{s}$ and $y\mid q+1$; $\frac{x}{s}\mid q-1$ and $q\mid y$; $\frac{x}{s}\mid q+1$ and $q\mid y$.

For the first case, $q=\frac{x}{s}$ and $\frac{q+ 1}{y}$ is an odd integer. This together with $\gcd(d, q(q^2-1))=1$ implies that $\gcd(y,xy-2x-2y)=\gcd(y,2)\leq 2$. Observing the equation $(\dag)$, we deduce that $y=q-1$ and there exist positive integers $j$ and $k$ such that $j+k=4$, $\frac{q+1}{2}=p^j$ and $xy-2x-2y=2p^k$. By the equality $\frac{q+1}{2}=p^j$, we have $q=2p^j-1$. Substituting $q=2p^j-1, x=sq=s(2p^j-1)$ and $y=q-1=2p^j-2$ in the equation $xy-2x-2y=2p^k$, we get $s=\frac{2p^k+4p^j-4}{(2p^j-1)(2p^j-4)}$. Similarly, for the remaining three cases, we obtain $s=\frac{2p^k+4p^j+4}{2p^j(2p^j+1)}$, $s=\frac{2p^k+4p^j-2}{(2p^j-2)(2p^j-3)}$ and $s=\frac{2p^k+4p^j+2}{(2p^j-1)(2p^j+2)}$, respectively. By Lemma \ref{int}, we see that $s$ cannot be an integer in all these situations. Thus case (2) cannot occur. This is a final contradiction, the claim follows. \qed
\vskip 2mm
Now we are able to exclude the case that $G$ is an insolvable almost Sylow-cyclic group. This is the following lemma.

\begin{lem}\label{excep}
Assume condition $(\star\star)$. Then $G$ cannot be an insolvable almost Sylow-cyclic group.
\end{lem}
\demo We keep the notation as in Lemma \ref{notdi} and write $\MM=\MM(G;r,t,\ell)$. Suppose for the contrary that $G$ is an insolvable almost Sylow-cyclic group. Then it follows from Lemma \ref{notdi} that $G$ is isomorphic to either $\PSL(2,p)$ or $\mathbb{Z}_d\rtimes_\phi \PGL(2,p)$ and $p\nmid xy$. Thus the equation $(\dag)$ in Lemma \ref{notdi} holds. In what follows, we consider these two cases separately.
\vskip 2mm
{\it Case $(1)$: $G\cong \PSL(2,p)$}
\vskip 2mm
In this case, we have $|G|=|\PSL(2,p)|=\frac{p(p^2-1)}{2}$. But on the other hand, we have $|G|= \frac{4xy}{xy-2x-2y} p^4>p^4>\frac{p(p^2-1)}{2}$, a contradiction.
\vskip 2mm
{\it Case $(2)$: $G\cong \mathbb{Z}_d\rtimes_\phi\PGL(2,p)$}
\vskip 2mm
In this case, $d>1$ and $\gcd(d,p(p^2-1))=1$. Without loss of generality, we may assume that $d\mid x$. Thus the subgroup $\mathbb{Z}_d\rtimes\PSL(2,p)$ contains an element of order $x$. In particular, $\frac{x}{d}\mid \frac{p\pm1}{2}$. Observing the equation $(\dag)$, we find that $\frac{p^2-1}{4}\mid \frac{x}{d}\cdot y$ is even and therefore $(\frac{x}{d},y)=(\frac{p-1}{2},p+1)$ or $(\frac{p+1}{2},p-1)$ and $xy-2x-2y=2p^3$. Solving the equation $xy-2x-2y=2p^3$, we obtain $d=\frac{4p^3+4p+4}{(p-1)^2}$ and $d=\frac{4p^3+4p-4}{(p+1)(p-3)}$ according to the two cases $\frac{x}{d}=\frac{p-1}{2}$ and $\frac{x}{d}=\frac{p+1}{2}$, respectively. Since $d$ is an integer, it follows that either $(p-1)^2\mid 4p^3+4p+4$ or $(p+1)(p-3)\mid  4p^3+4p-4$. We deduce that $p-1\mid 12$ when $(p-1)^2\mid 4p^3+4p+4$ and $p+1\mid 12$ when $(p+1)(p-3)\mid  4p^3+4p-4$. Since we are assuming $p\geq 5$, we derive that $p=5, 7, 13$ when  $p-1\mid 12$ and $p=5, 11$ when $p+1\mid 12$. Clearly,  we have only two cases for when $d$ is an integer: $p=5$ and $d=43$; $p=7$ and $d=39$. For the case $p=7$ and $d=39$, we have $\gcd(39, 7(7^2-1))=3>1$, a contradiction to $\gcd(d, q(q^2-1))=1$. Thus we are left with the case $p=5$ and $d=43$. In this case, it is not hard to see that $\{x, y\}=\{4, 129\}$. Let $\overline{\MM}$ be the quotient map of $\MM$ induced by the normal subgroup of order $d$. Then $\overline{\MM}$ has type $\{3, 4\}$ and $\Aut(\overline{\MM})\cong \PGL(2,5)$. However, the order of the automorphism group of a regular map with type $\{3, 4\}$ is at most $48$. This is a contradiction to the fact that $48<120=|\PGL(2,5)|$. \qed 

\vskip 2mm

Assume condition $(\star\star)$. The following lemma gives an upper bound of the index of $O_p(G)$ in $P$. In particular, it follows from Proposition \ref{alm} and Theorem \ref{divide} that $G/O_p(G)$ is an almost Sylow-cyclic group.
\vskip 2mm
\begin{lem}\label{nontri}
Assume condition $(\star\star)$ and let $P\in \Syl_p(G)$. Then $|P:O_p(G)|\leq p$.
\end{lem}
\demo We abbreviate $\MM(G;r,t,\ell)$ by $\MM$ and let $\{x,y\}$ be the type of $\MM$. The statement is trivial if $P$ is trivial. Thus we may assume that $P$ is nontrivial, that is, $p\mid |G|$. Assume that $P>1$ is cyclic. Then it follows from Proposition \ref{alm} that $G$ is almost Sylow-cyclic. According to Lemma \ref{excep}, we see that $G$ is solvable. But by Proposition \ref{solvable}, we know that $p\nmid |G|$, a contradiction. Thus we may assume that $P$ is non-cyclic. In particular, we have $|P|=p^2$ or $p^3$. By Proposition \ref{normal-sub}, we have $|O_p(G)|\geq p$. If  $|P|=p^2$, then $|P:O_p(G)|\leq \frac{p^2}{p}=p$, the statement follows. Thus we only need to consider the case $|P|=p^3$. 

Assume $|P|=p^3$. All we need to show is that $O_p(G)\not\cong \mathbb{Z}_p$. Suppose for the contrary that $O_p(G)\cong \mathbb{Z}_p$. Then $G/O_p(G)$ has nontrivial Sylow $p$-subgroup. In particular, it follows from Proposition \ref{p^3} that the Euler characteristic of the quotient map of $\MM$ induced by $O_p(G)$ is never $-p^3$. This implies that $O_p(G)$ is contained in either $\lg rt\rg$ or $\lg r\ell \rg$. It follows from $|P|=p^3$ and formula $(*)$ that $p\mid x$ and $p\mid y$. Let $L$ denote the largest normal subgroup of $G$ of odd order. We claim that $p^2\nmid |L|$. To see it, note that it follows from $O_p(G)=O_p(L)$ that $F(L)\leq L$ is a cyclic group. Note that any Sylow $l$-subgroup of $G$ with $l\neq p$ an odd prime is cyclic. Thus it follows from
$G=\lg rt, t\ell\rg=\lg r\ell, t\ell\rg$ and Proposition \ref{nc} that $L\leq C_G(F(L))$. But on the other hand, the solvability of $L$ implies that $C_L(F(L))\leq F(L)\leq L$. Hence $F(L)=L$ is cyclic. This implies that $O_p(G)$ is a Sylow $p$-subgroup of $L$. Therefore, $p^2\mid |G/L|$. According to Proposition \ref{dihedral}, $G/L$ is isomorphic to an almost simple group with socle $\PSL(2,d)$. Clearly, one can see that $p^2\leq d$. Thus $|G|=|L||G/L|\geq p\cdot \frac{d(d^2-1)}{2}\geq p\cdot \frac{p^2(p^4-1)}{2}$. But on the other hand, we have $|G|\leq 84p^4$. Hence $p\cdot \frac{p^2(p^4-1)}{2}\leq 84p^4$. Under our assumption, this inequality holds if and only if $p=5$. With a similar argument as in Proposition \ref{normal-sub}, we deduce that $$\frac{1}{x}+\frac{1}{y}\geq \frac{1}{3}+\frac{1}{13}=\frac{16}{39}.$$But then $$\frac{5^3(5^4-1)}{2}\leq |G|=\frac{4}{1-2(\frac{1}{x}+\frac{1}{y})}\cdot 5^4\leq 4\cdot \frac{39}{7}\cdot 5^4,$$ a contradiction. \qed

\vskip 2mm
The following lemma will be useful in our later proof. 
\begin{lem}\label{noi}
Let $p\geq 5$ be a prime. Then the following hold:

$(1)$ The equation $pmn-km-kn=dp^2$ with $d,k\in \{1,2\}$ has positive integer solutions for $(m,n)$ if and only if 

$(a)$ $d=1, p=5$, and either $k=1$ and $\{m,n\}=\{2,3\}$, or $k=2$ and $\{m,n\}=\{1,9\}$;
\vskip 2mm
$(b)$ $d=2, p=7$, $k=2$ and $\{m,n\}=\{1,20\}$.
\vskip 2mm
$(2)$ The equation $pmn-m-n=p^3$ has no positive integer solutions for $(m,n)$ such that $2\nmid mn$.
\end{lem}
\demo For convenience, we abbreviate the equation $pmn-km-kn=dp^2$ by Equation $(\ddag)$. Suppose $m$ and $n$ are positive integers solutions to Equation $(\ddag)$. Without loss of generality, we may assume that $m\leq n$. We first consider the case $d=1$. Suppose $(m-2)(n-2)\leq 0$. Then either $m=1$ or $m=2$. Suppose  $m=1$. Then $pn-k-kn=p^2$. Thus $n=p+k+\frac{k^2+k}{p-k}$. In particular, $\frac{k^2+k}{p-k}$ is an integer. This together with the assumption $p\geq 5$ yields that $p=5$ and $k=2$. Hence $n=9$ in this case. Suppose $m=2$. Then we see from the parity of the both sides of the Equation $(\ddag)$ that $k=1$. Thus from the Equation $(\ddag)$ we obtain $4n=2p+1+\frac{9}{2p-1}$. In particular, $\frac{9}{2p-1}$ is an integer. This together with the assumption $p\geq 5$ yields that $p=5$. In this case, $n=3$. Assume $(m-2)(n-2)\geq 6$, equivalently, $2(m+n)\le mn-2$. Thus we have $pmn-(mn-2)\le pmn-k(m+n)=p^2$. Hence $mn<p+1-\frac{1}{p-1}<p+1$. On the other hand, we deduce from the Equation $(\ddag)$ that $p^2<p^2+k(m+n)=pmn$. This implies that $mn>p$. Therefore, $mn$ cannot be an integer, a contradiction. Now, we are left with the case $1\leq (m-2)(n-2)\leq 5$. Then either $m=3$ and $n\in \{3, 4, 5, 6, 7\}$, or $m=n=4$. With a straightforward case-by-case analysis, we find that $m$ and $n$ are never solutions to the Equation $(\ddag)$ for these cases. Similarly, we deduce that $p=7$, $k=2$ and $\{m,n\}=\{1,20\}$ for $d=2$.

Suppose for the contrary that the equation $pmn-m-n=p^3$ has positive integer solutions for $(m,n)$ such that $2\nmid mn$. Then $m$ and $n$ are odd. In particular, $(m-1)(n-1)\equiv 0({\rm mod}~4)$. Thus $$(p-1)mn\equiv (p-1)mn+(m-1)(n-1) \equiv p^3+1({\rm mod}~4).$$ Since we are assuming $p\geq 5$, we have either $p\equiv 1({\rm mod}~4)$ or $p\equiv 3({\rm mod}~4)$. Assume $p\equiv 1({\rm mod}~4)$. Then $(p-1)mn\equiv 0({\rm mod}~4)$ while $p^3+1\equiv 2({\rm mod}~4)$, a contradiction. Assume $p\equiv 3({\rm mod}~4)$. Then $(p-1)mn\equiv 2({\rm mod}~4)$ while $p^3+1\equiv 0({\rm mod}~4)$, a contradiction. \qed
\vskip 2mm
The following lemma will be useful in our later proof.
\begin{lem}\label{common}
Assume condition $(\star\star)$ and suppose that $G$ is insolvable. Then $\gcd(|rt|_{p'}, |r\ell|_{p'})>1$.
\end{lem}
\demo Let $x$ and $y$ denote $|rt|$ and $|r\ell|$, respectively, and let $x'$ and $y'$ denote $x_{p'}$ and $y_{p'}$, respectively. Suppose for the contrary that $\gcd(x', y')=1$. Then $2$ cannot divide $x'$ and $y'$ simultaneously. This implies that at least one of $x'$ and $y'$ is odd. Without loss of generality, we may assume that $y'$ is odd. 

Let $P\in \Syl_p(G)$, and let $p^k$ and $p^s$ be the order of $P$ and the $p$-part of $\gcd(x,y)$, respectively. Then $k\leq 3$ and it follows from formula $(*)$ that $x=p^s\cdot x'$ and $y=p^s\cdot y'$. By Lemma \ref{excep}, we know that $P$ is non-cyclic. In particular, $k\geq 2$ and neither $\lg rt\rg$ nor $\lg r\ell\rg$ contains a Sylow $p$-subgroup of $G$. This implies that $p^s<p^k$, and therefore $s\leq k-1\leq 2$. Since $|G:\lg rt\rg|=\frac{|G|}{2x}=\frac{2y'p^4}{p^sx' y'-2x'-2y'}$ is a positive integer the $p$-part of which is exactly $p^{k-s}$ and $\gcd(y',p^sx' y'-2x'-2y')=\gcd(y',2x')=1$, we obtain one of the following:

$(1)$ $x'$ is odd and $p^sx' y'-2x'-2y'=p^{4+s-k}$;
\vskip 2mm
$(2)$ $x'$ is even and $p^sx' y'-2x'-2y'=2p^{4+s-k}$. 

\noindent Thus $|G|=\frac{4xyp^{k-2s}}{\gcd(2,x')}$. In particular, $|G|\leq 4xyp^{k-2s}$.

We show that $s\neq 0$. Suppose for the contrary that $s=0$. Then $p\nmid xy$. Thus the quotient map of $\MM$ induced by $O_p(G)$ is of type $\{x,y\}$ and of Euler characteristic $-p^j=\frac{-p^4}{|O_p(G)|}$. Since $k\geq 2$, it follows from Proposition \ref{normal-sub} that $O_p(G)>1$. This together with $O_p(G)\leq P$ implies that $1\leq j<4$. Assume $j=1$. Then $|O_p(G)|=p^3$. Since $G/O_p(G)$ is insolvable, it follows from Proposition \ref{cla} that $G/O_p(G)$ has nontrivial Sylow $p$-subgroup. This yields that $P>O_p(G)$ and therefore $|P|\geq p^4$, a contradiction to our assumption. Assume $j=2$. Then by Proposition \ref{Eule} that $p=7, \{x, y\}=\{3, 13\}$ and $G/O_p(G)\cong \PSL(2,13)$. According to Proposition \ref{nc}, we deduce that $O_p(G)\leq Z(G)$. Therefore, $O_P(G)\leq Z(P)$. As $O_P(G)\leq Z(P)$ and $P/O_p(G)$ is cyclic, it follows that $P$ is abelian. By Proposition \ref{is}, we see that $O_p(G)\cap G'=1$. In particular, $[O_p(G), G']=1$ and $G'<G$. Clearly, $G'>1$ as $G$ is insolvable. Thus we obtain from $1<(O_p(G)G')/O_p(G)\unlhd G/O_p(G)$ that $G=O_p(G)\times G'$. But then $rt, r\ell\in G'$ and therefore $G=\lg rt, r\ell\rg\leq G'<G$, a contradiction. Assume $j=3$. Then by Proposition \ref{p^3} we derive that $G/O_p(G)$ is solvable. This implies that $G$ is solvable, a contradiction to our assumption. Hence we have $s\geq 1$ and therefore $(s, k)=(1,2), (1,3)$ or $(2,3)$.

We claim that $(s, k)=(1,2)$. If this were proved, then $|G|=\frac{4xy}{\gcd(2,x')}$. Suppose for the contrary that $(s, k)=(1,3)$ or $(2,3)$. Assume $(s, k)=(2,3)$. Then $G\leq \frac{4xy}{p}<xy$ as $p\geq 5$. But on the other hand, it follows from $\lg rt\rg\cap \lg r\ell\rg=1$ and $\lg rt\rg\lg r\ell\rg\subseteq G$ that $|G|\geq xy$. This is a contradiction. Assume $(s, k)=(1,3)$. Then by Lemma \ref{noi} we have either $p=5$ and $\{x',y'\}=\{1,9\}$, or $p=7$ and $\{x',y'\}=\{1,20\}$. Thus either $p=5,\{x,y\}=\{5,45\}$ and $|G|=2^2\cdot 3^2\cdot 5^3$, or $p=7,\{x,y\}=\{7,140\}$ and $|G|=5\cdot 2^3\cdot 7^3$. Since we are assuming that $G$ is insolvable, it follows from Proposition \ref{dihedral} that $3\mid |G|$. This implies that $|G|=2^2\cdot 3^2\cdot 5^3$. By comparing the order of $|G|$ with that of almost simple groups $N$ with socle $\PSL(2, f)$, we deduce that $G/O_{2'}(G)\cong \PSL(2,5)$ and $|O_{2'}(G)|=3\cdot 5^2$. According to Lemma \ref{nontri}, we have $|O_5(G)|\geq 5^2$. This together with $G/O_{2'}(G)\cong \PSL(2,5)$ implies that $G/O_5(G)\cong \mathbb{Z}_3\times \PSL(2,5)$. On the other hand, it follows from Proposition \ref{alm} that $G/O_5(G)$ is almost Sylow-cyclic. However,  $G/O_5(G)\cong \mathbb{Z}_3\times \PSL(2,5)$ is not almost Sylow-cyclic, a contradiction. The claim follows.

We show that $x'$ is odd. Set $\og:=G/O_{2'}(G)$ and $\widehat{G}:=G/O_p(G)$. Then there is the canonical surjective homomorphism $\varphi:\widehat{G}\twoheadrightarrow \og, \hat{g}\mapsto \overline{g}$ for all $g\in G$. Note that $\langle rt\rangle\cap \langle r\ell\rangle=1$. Then we obtain $|\langle \olr\olt\rangle \langle \olr\oll\rangle|\geq \frac{\gcd(2,x')\cdot |\og|}{4}$ and $|\langle \widehat{rt}\rangle \langle \widehat{r\ell}\rangle|\geq \frac{\gcd(2,x')\cdot |\widehat{G}|}{4}$. Recall that $G$ has dihedral Sylow-$2$ subgroup by Proposition \ref{alm}. According to Proposition \ref{dihedral}, $\og$ is isomorphic to either either an almost simple group with socle $\PSL(2, q)$ where $q$ is odd, or $A_7$, or a Sylow $2$-subgroup of $G$. It is straightforward to check that $A_7$ cannot be generated by three involutions, two of which commute. Thus we can exclude the case $\og\cong A_7$. Since we are assuming the solvability of $G$, it follows easily that $\og$  cannot be a $2$-group. Therefore, $\og$ is isomorphic to an almost simple group with socle $\PSL(2, q)$, where $q\geq 5$ is a prime power. On the other hand, $\widehat{G}$ is an almost Sylow-cyclic insolvable group by Lemma \ref{nontri}. Then $\og$ is also almost Sylow-cyclic as it is a homomorphic image of $\widehat{G}$. This implies that $q$ is a prime and therefore $\PSL(2,q)\leq \og\leq \PGL(2,q)$. Let $\PSL(2,q)\leq N\leq \PGL(2,q)$ and let $R, S\leq N$ be two arbitrary cyclic subgroups. Using the well-known result on orders of elements in $\PGL(2,q)$, we have the following fact:

$$|R||S|<\frac{|N|}{4} ~\mbox{for}~ q\geq 7~\mbox{and}~|R||S|<\frac{|N|}{2}~\mbox{for}~q=5.$$This fact together with $|\langle \widehat{rt}\rangle \langle \widehat{r\ell}\rangle|\geq \frac{\gcd(2,x')\cdot |\widehat{G}|}{4}$ implies that $q=5$ and $x'$ is odd.

Since $x'$ is odd and $(s,k)=(1,2)$, we have $px' y'-2x'-2y'=p^3$ by $(1)$. But by Lemma \ref{noi}(2), this equation has no positive integer solutions for $(x',y')$ such that $2\nmid x'y'$. Thus we arrive a final contradiction.  \qed
\vskip 2mm
Now, with all the machineries needed, we are in a position to prove Theorem \ref{so-norm}.
\vskip 2mm
\noindent{\bf Proof of Theorem \ref{so-norm}.} Let $\{x,y\}$ be the type of $\MM$ and $x', y'$ be the $p'$-part of $x$ and $y$, respectively. By proposition \ref{alm} and Lemma \ref{nontri}, the quotient map of $\MM$ induced by $O_p(G)$, say $\overline{\MM}$, is a non-orientable regular map with almost Sylow-cyclic automorphism group $\og:=G/O_p(G)$. Let $\{x_1,y_1\}$ be the type of $\overline{\MM}$. Then it is clear that $x'\mid x_1$ and $y'\mid y_1$. Note that $\og\cong \PSL(2,q)$ or $\mathbb{Z}_d\rtimes \PGL(2,q)$ for some prime $q$ and integer $d$ such that $\gcd(d, q(q^2-1))=1$. Without loss of generality, we may assume that $d\mid x_1$ in the latter case. We only need to show that $G$ is solvable. If this were proved, then $\Aut(\overline{\MM})=\og$ would be a solvable almost Sylow-cyclic group, and therefore $\og$ has a normal Sylow $p$-subgroup by Proposition \ref{alm}, whence $P=O_p(G)\unlhd G$. Suppose for the contrary that $G$ is insolvable. In what follows, we prove $q=p$ and $O_p(G)\cong \mathbb{Z}^2_p$. 

We first show that $q=p$. Suppose for the contrary that $q\neq p$. Since any odd prime divisor of $|\og|$ except for $p$ divides $x_1 y_1$, it follows from $3q\mid |\og|$ that $3q$ divides $x_1y_1$. This together with $\og\cong \PSL(2,q)$ or $\mathbb{Z}_d\rtimes \PGL(2,q)$ implies that $q\in \{\frac{x_1}{\alpha},y_1\}$ and $3$ divides the element in $\{\frac{x_1}{\alpha},y_1\}\setminus \{q\}$, where $\alpha=1$ if $\og \cong \PSL(2,q)$ and $\alpha=d$ if $\og\cong \mathbb{Z}_d\rtimes \PGL(2,q)$. Since the element in $\{\frac{x_1}{\alpha},y_1\}\setminus \{q\}$ is divided by $3$, it follows that this element divides $(q-1)(q+1)$. Thus $\gcd(x', y')\leq \gcd(x_1,y_1)\leq \gcd(q, (q+1)\cdot (q-1))=1$. This is a contradiction with Lemma \ref{common}, as desired.

Next we show that $O_p(G)\cong \mathbb{Z}^2_p$. By Proposition \ref{non-alm}, $\og\cong \PSL(2,q)$ or $\mathbb{Z}_d\rtimes_\phi \PGL(2,q)$ (here we may assume without loss of generality $d\mid x_1$), where $q\geq 5$ is a prime, $d$ is an integer coprime to $q(q^2-1)$ and $\phi: \PGL(2,q)\twoheadrightarrow \langle \sigma_{\mathbb{Z}_d}\rangle$ is a homomorphism. In particular, $1, 2\notin \{x_1, y_1\}$ as $\og$ is never a dihedral group. Assume that $P$ is normal. Then $\{x_1,y_1\}=\{x',y'\}$ and any odd prime divisor of $|\og|$ divides $x_1 y_1$. Note that $\og=\langle \olr\olt,\olr\oll\rangle$. Applying Lemma \ref{cop}, we conclude that $\gcd(x',y')=1$, a contradiction with Lemma \ref{common}. Thus $P>O_p(G)$. Since $|P|\leq p^3$ by our assumption, it follows that $|O_p(G)|\leq p^2$. We claim that $p\mid xy$. Suppose for the contrary $p\nmid xy$. Then $\{x_1,y_1\}=\{x,y\}$ and the quotient map $\overline{\MM}$ is a regular map with $\chi_{\overline{\MM}}=-\frac{p^4}{|O_p(G)|}$. Assume $|O_p(G)|=p^2$. Then $\chi_{\overline{\MM}}=-p^2$. By Proposition \ref{Eule}, we deduce that $p=7$, $\og \cong \PSL(2,13)$ and $\{x_1,y_1\}=\{3,13\}$. Applying Proposition \ref{nc}, we derive that $G=O_p(G)\times H$ for some subgroup $H\cong \PSL(2,13)$. But then any element of order $3$ or $13$ lies in $H$, a contradiction to the fact that $G=\lg rt, r\ell\rg$ is generated by an element of order $3$ and an element of order $13$.  Assume $|O_p(G)|=p$. Then $\chi_{\overline{\MM}}=-p^3$. By  Proposition \ref{p^3} we deduce that $G/O_p(G)$ is solvable and therefore $G$ is solvable, a contradiction to our assumption. The claim follows. This together with the formula $(*)$ and $|P|\leq p^3$ implies that $x_p=y_p\geq p$. Assume that $O_p(G)$ intersects trivially with at least one of $\langle rt\rangle$ and $\langle r\ell\rangle$. Then $p\mid x_1y_1$. Therefore, any odd prime divisors of $|\og|$ divides $x_1y_1$. By Proposition \ref{cop}, we derive that $\gcd(x_1,y_1)=1$. This implies that $\gcd(x',y')=1$, a contradiction with Lemma \ref{common}. Thus $O_p(G)$ intersects nontrivially with both $\langle rt\rangle$ and $\langle r\ell\rangle$. Since $\lg rt\rg \cap \lg r\ell\rg=1$ and $|O_p(G)|\leq p^2$, it follows that $O_p(G)\cong \mathbb{Z}^2_p$, as desired. 

Since $q=p$, it follows that $|\og|_{p'}$ divides $4x_1 y_1$. This together with $\gcd(x_1,y_1)\geq \gcd(x', y')>1$ yields that there are the following three situations:

(a) $G/O_p(G)\cong \PSL(2,p)$, both $x_1$ and $y_1$ divide $\frac{p-1}{2}$ (resp. $\frac{p+1}{2}$) and $\frac{p+1}{2}$ (resp. $\frac{p-1}{2}$) is a divisor of $8$

(b) $G/O_p(G)\cong \mathbb{Z}_d\rtimes_\phi \PGL(2,p)$, both $\frac{x_1}{d}$ and $y_1$ divide $p-1$ (resp. $p+1$) and $p+1$ (resp. $p-1$) is a divisor of $8$

(c) $G/O_p(G)\cong \mathbb{Z}_d\rtimes_\phi \PGL(2,p)$, $\frac{x_1}{d}=\frac{p+1}{k}$, $y_1=\frac{p-1}{l}$ or $\frac{x_1}{d}=\frac{p-1}{k}$, $y_1=\frac{p+1}{l}$, where $kl\mid 4$.

Assume case (a). Then we have $p=7$ and $\{x_1,y_1\}=\{3,3\}$, or $p=17$ and $x_1,y_1\in \{3,9\}$. Assume case (b). Then we have $p=5, 7$ and $\frac{x_1}{d},y_1\in \{3,6\}$, or $p=17$ and $\frac{x_1}{d},y_1\in \{3,6,9,18\}$. Substituting all this possibilities for $x=px_1, y=py_1$ and $p$ in the formula $(*)$, we find a contradiction. Assume case (c). Then by formula $(*)$ we obtain $\frac{dp(p^2-1)}{kl}-\frac{2d(p\pm 1)}{k}-\frac{2(p\mp 1)}{l}=p^2$. Assume $d>1$. Then it follows from $\gcd(d,p(p^2-1))=1$ that $d\geq 5$. Thus $$\small{p^2=\frac{dp(p^2-1)}{kl}-\frac{2d(p\pm 1)}{k}-\frac{2(p\mp 1)}{l}\geq \frac{dp(p^2-1)}{kl}-\frac{2d(p\pm 1)}{k}-\frac{2d(p\mp 1)}{k}=\frac{dp(p^2-1-4l)}{kl}.}$$ Since we are assuming $p\geq 5$, it follows that $p^2-1-4l\geq p^2-9\geq p$. This implies that $\frac{dp(p^2-1-4l)}{kl}\geq \frac{dp^2}{4}>p^2$, a contradiction. Then $d=1$. Considering the four cases $(k,l)=(1,1),(1,2),(2,1)$ and $(2,2)$ separately, we find that $(k,l)=(1,1), p=5$ and $\{x_1, y_1\}=\{4,6\}$. In this case, $\MM$ is of type $\{20, 30\}$ and $G$ is an extension of $\PGL(2,5)$ by $\mathbb{Z}_5\times \mathbb{Z}_5$. But this cannot happen by Lemma \ref{nonexi}, a final contradiction. \qed
\begin{lem}\label{nonexi}
The automorphism group of a regular map is never isomorphic to an extension of $\PGL(2,5)$ by $\mathbb{Z}_5\times \mathbb{Z}_5$. 
\end{lem}
\demo Let $G$ be an extension of $\PGL(2,5)$ by $\mathbb{Z}_5\times \mathbb{Z}_5$. We first show that $G$ has a subgroup, say $H$, such that $G=O_5(G)\rtimes H$ and $H\cong \PGL(2,5)$. To see it, let $K\leq G$ be the subgroup containing $O_5(G)$ such that $K/O_5(G)\cong \PSL(2,5)$. Then by Proposition \ref{nc} we derive that $O_5(G)\leq Z(K)$ as $\PSL(2,5)$ is not isomorphic to a subgroup of $\GL(2,5)$. By a similar argument as for $s\neq 0$ in the proof of Lemma \ref{common}, we see that $K=O_5(G)\times K'$. Clearly, $K'\cong \PSL(2,5)$. As $K'\char K\unlhd G$, it follows that $K'\unlhd G$. Since $K<G=\lg r, t,\ell\rg$, it follows that at least one of $r, t$ and $\ell$ lies out of $K$, and denote it by $z$. Then $K'\lg z\rg\cap O_5(G)=1$ and $K'\lg z\rg O_5(G)=K\lg z\rg=G$. This together with $G/O_5(G)\cong \PGL(2,5)$ implies that $K'\lg z\rg\cong \PGL(2,5)$. Thus $H:=K'\lg z\rg$ is our desired subgroup. Clearly, $[z, O_5(G)]\neq 1$, otherwise $H\unlhd G$ is a proper subgroup containing all involutions of $G$, which is a contradiction to $G=\lg r, t,\ell\rg$. 

Suppose for the contrary that $\MM=\MM(G;r,t,\ell)$ is a regular map. As $2\mid |G:H|=|O_5(G)|$, $H$ contains a Sylow $2$-subgroup of $G$. Thus by Sylow's theorems and using conjugation if necessary, we may assume that $\lg t,\ell\rg\leq H$. Since $G=O_5(G)\rtimes H$, we can write $r=gh$ for unique $g\in O_5(G)$ and $h\in H$. Clearly, it follows from $r^2=(gh)^2=1$ that $g^h=g^{-1}$. Clearly, $g\neq 1$, otherwise $H<G=\lg r, t, \ell\rg =\lg h, t, \ell\rg \leq H$, a contradiction. This together with $[K', O_5(G)]=1$ implies that $h\notin K'$. Thus $H=K'\lg h\rg$ and therefore $H$ normalizes $\lg g\rg$. In particular, $h, t$ and $\ell$ all normalize $\lg g\rg$. Hence $$G=\lg r, t, \ell\rg=\lg gh,t, \ell\rg\leq \lg g\rg \lg h, t,\ell\rg \leq \lg g\rg H\leq G.$$This yields that $G=\lg g\rg H$. However, comparing the order of $G$ and $\lg g\rg H$, we find a contradiction. \qed
\begin{rem}
As pointed out before Proposition \ref{cla}, \cite[Theorem 2.2(3)]{ARJ} excludes the case-a regular map of type $\{4, 6\}$ with Euler characteristic $-5$ and automorphism group $\PGL(2,5)$. Consequently, the case $p=5$ and $G/O_5(G)\cong \PGL(2,5)$ was not considered in Step 2 of the proof of \cite[Theorem 3.1]{TL}. But fortunately, it follows from Lemma \ref{nonexi} that this case is impossible.
\end{rem}

\subsection{The analysis of $P$}
In this subsection, we mainly prove:
\begin{theorem}\label{cyclic}
Assume condition $(\star\star)$ and let $P\in \Syl_p(G)$. Then $P$ is trivial.
\end{theorem}
\demo Let $\{x,y\}$ be the type of $\MM=\MM(G;r,t,\ell)$ and let $p^d$ be the order of $P$. Then $d<4$ by our assumption. We claim that $G$ is almost Sylow-cyclic. If claim were proved, then the statement would follow immediately by Proposition \ref{solvable}. Suppose for the contrary that $G$ is not almost Sylow-cyclic. Then Proposition \ref{alm} implies that $P$ cannot be cyclic. In particular, $d\geq 2$. Note that $P\unlhd G$ by Theorem \ref{so-norm}. Let $\overline{\MM}$ be the quotient map of $\MM$ induced by $P$ and let $x':=x_{p'}$ and $y':=y_{p'}$. Then $\overline{\MM}$ is of type $\{x', y'\}$ and it follows from Proposition \ref{alm} and Theorem \ref{so-norm} that $\Aut(\overline{\MM})=G/P$ is a solvable almost Sylow-cyclic group. Applying Proposition \ref{solvable}, one of the following holds
\begin{enumerate}
\item[\rm (1)] $\{x',y'\}=\{2,n\}$, $t$ even, $G/P\cong \mathbb{D}_{2n}$;
\item[\rm (2)] $\{x',y'\}=\{2a,2b\}$, $G/P\cong G(a,b)$, where $2\nmid ab$, $1\neq a<b$ and $\gcd(a,b)=1$;
\item[\rm (3)] $\{x',y'\}=\{4,j\}$, $G/P\cong G(j)$, where $j\equiv 3({\rm mod}~6)$.
\end{enumerate}
We first show that $p\nmid xy$. To see it, we use contradiction and assume $p\mid xy$. Then it follows from formula $(*)$ and the assumption $|P|<p^4$ that $x_p=y_p>1$. Let $R\in \Syl_p(\lg rt \rg)$ and $Q\in \Syl_p(\lg r\ell \rg)$. Then $R\cap Q\leq \lg rt \rg \cap \lg r\ell\rg=1 $ and $RQ\subseteq P$. In particular, it follows that $|R||Q|\leq |P|<p^4$. This yields that $x_p=y_p=p$ and therefore $2\leq d<4$. For the case (1), $\{x,y\}=\{2p,pn\}, |G|=2p^dn$; for the case (2), $\{x,y\}=\{2pa,2pb\}, |G|=4p^dab$; for the case (3), $\{x,y\}=\{4p,jp\}, |G|=8p^dj$. Substituting $x,y$ and $|G|$ into the formula $(*)$, we find that $n=2(\frac{p^{5-d}-1}{p-1}+\frac{2}{p-1})$ in the case (1), $pab-a-b=p^{4-d}$ in the case (2) and $2pj-j=p^{5-d}+4$ in the case (3). The first case cannot happen as $\frac{2}{p-1}$ is never an integer for $p\geq 5$ and the second case is excluded by Lemma \ref{noi}. For the third case, as the left hand side of the equality $2pj-j=p^{5-d}+4$ is divided by $3$, we have $d=2$ and 
$p\equiv 2({\rm mod}~3)$. But then the left hand side of the equality $2pj-j=p^{5-d}+4$ is divided by $9$ while the right hand side is not, a contradiction. Hence $p\nmid xy$. This implies that $\overline{\MM}$ is a regular map with Euler characteristic $-p^{4-d}$. Since $G$ is solvable, it follows from Proposition \ref{Eule} that $d=3$ and $\chi_{\overline{\MM}}=-p$. Thus $|P|=p^3$. Suppose that $P$ is not elementary abelian. Then $P$ always has a characteristic subgroup, say $K$, of order $p$. Consider the quotient map, say $\widetilde{\MM}$ of $\MM$, induced by $K$. Then $\chi_{\widetilde{\MM}}=-p^3$. By Proposition \cite[Theorem 1.1]{TL}, we derive $p\nmid |G/K|$, a contradiction to $K<P$. Hence $P$ is elementary abelian of rank $3$ and only the case (2) and the case (3) can occur. 

Assume case (2). Then we have $ab-a-b=p$, $G=P\rtimes H$ for some subgroup $H\cong G(a,b)$ and $\MM$ is of type $\{2a,2b\}$. Consider the conjugation of $H$ on $P$. We claim that $P$ has no nontrivial $H$-invariant subgroup. To see it, suppose that $P$ has a nontrivial $H$-invariant subgroup, say $T$. Then $T\unlhd G$ and so the quotient map of $\MM$ induced by $T$ is of Euler characteristic $-p^2$ or $-p^3$. But by Proposition \ref{Eule} and \cite[Theorem 1.1]{TL}, we know that $p\nmid |G/T|$. This implies $T=P$, a contradiction to $T<P$. Let $A$ and $B$ be the subgroup of $H$ of order $a$ and $b$, respectively. We claim that neither $A$ nor $B$ acts trivially on $P$. To see it, we may suppose $A$ acts trivially on $P$ (the analysis for $B$ acting trivially on $P$ is similar). Then $A\unlhd G$ and the quotient map of $\MM$ induced by $A$ is of type $\{2,2b\}$. This implies that $$4\cdot 5^3\cdot ab\leq 4p^3ab=|G|=|A||G/A|\leq a\cdot 8b,$$ which is clearly a contradiction. According to Proposition \ref{nc}, we know that $H/C_H(P)\cong G(a',b')$ is isomorphic to an irreducible subgroup of $\GL(3,p)$, where $a'>1$ and $b'>1$ are divisors of $a$ and $b$, respectively. Due to the well-known Clifford theory in group representation, we deduce that the normal subgroup $L$ of $G(a',b')$ of order $a'b'$ is irreducible. Recall that each irreducible cyclic subgroup of a general linear group is a subgroup of some Singer cycle group (see for instance \cite[Theorems 2.32, 2.33]{MKS}). But by Proposition \ref{HU}, the normalizer of an irreducible cyclic subgroup of $\GL(3,p)$ is isomorphic to $\mathbb{Z}_{p^3-1}\rtimes \mathbb{Z}_3$. Thus $G(a',b')$ is isomorphic to a subgroup of $\mathbb{Z}_{p^3-1}\rtimes \mathbb{Z}_3$, which is clearly a contradiction. 

Assume case (3). Then we have $j-4=p$, $G=P\rtimes S$ for some subgroup $S\cong G(j)$ and $\MM$ is of type $\{4,j\}$. By the same reason as in the case (2), we know that $S$ acts irreducibly on $P$. Let $D$ be the normal subgroup of $S$ of order $4$. Suppose that $D$ acts trivially on $P$. Then $D\unlhd G$ and so the quotient map of $\MM$ induced by $D$ is of type $\{2,j\}$. Thus $$8\cdot 5^3\cdot j\leq 8p^3j=|G|=|D||G/D|\leq 4\cdot 4j,$$ which is clearly a contradiction. According to Proposition \ref{nc}, we get that $\os=S/C_S(P)\cong G(j')$ is isomorphic to an irreducible subgroup of $\GL(3,p)$, where $j'$ is a divisor of $j$. Without loss of generality, we may assume that $S/C_S(P)\leq GL(3,p)$.   Since $\gcd(p+4,p+1)=\gcd(p+4,3), \gcd(p+4,p-1)=\gcd(p+4,5)$ and $\gcd(p+4,p^2+p+1)=\gcd(p+4,13)$, we know from $|\GL(3,p)|=p^3(p-1)^3(p+1)(p^2+p+1)$ that $j'=3, 15$ or $39$. For the case $j'=3$ (note that $j'=3$ if $p=5$), the quotient map of $\MM$ induced by $C_S(P)$ is of type $\{3,4\}$. This implies that $$8\cdot 5^3\cdot j\leq 8p^3j=|G|=|C_S(P)||G/C_S(P)|\leq \frac{j}{3}\cdot 48,$$ which is clearly a contradiction. For the case $j'=39$, let $\oi\leq \os$ be the subgroup of order $13$. As $13$ divides $p^3-1$ but not $p^2-1$, we deduce that $\oi$ is an irreducible subgroup of $\GL(3, p)$. By Proposition \ref{HU}, the centralizer of  $\oi$ in $GL(3,p)$ is cyclic of order $p^3-1$, a contradiction to the fact that $\od$ centralizes $\oi$ but is not cyclic. Finally, we are left with the case $j'=15$ so that $p\equiv 1 ({\rm mod}~5)$. Write $\od=\lg a,b\rg$. Then we can find elements $u,v,w\in P$ such that $P=\lg u\rg \times \lg v\rg\times \lg w\rg$ decomposes as an $\od$-module. By the irreducibility of $P$ as an $\os$-module, we know that $\od$ fixes no nonidentity elements of $P$. Thus we may assume that $a$ and $b$ are represented by the diagonal matrices ${\rm diag}(1,-1,-1)$ and ${\rm diag}(-1,1,-1)$ with respect to the basis $u,v$ and $w$. Let $c\in \os$ be  an element of order $5$. Then $c$ commutes with $\lg a,b\rg$ and it follows that $\lg u\rg $, $\lg v\rg $ and $\lg w\rg $ are stabilized by $c$. Thus we may assume that $c$ is represented by a matrix ${\rm diag}(\lambda_1,\lambda_2,\lambda_3)$ with respect to the basis $u,v$ and $w$. Since $c\notin Z(\os)$, ${\rm diag}(\lambda_1,\lambda_2,\lambda_3)$ can not be a scalar matrix. Thus we may assume that $\lambda_1$ is different from $\lambda_2$ and $\lambda_3$. Let $d\in \os$ be an element of order $3$. Since $c$ and $d$ commutes, it follows from $\lambda_1\notin \{\lambda_2,\lambda_2\}$ that $d$ stabilizes $\lg u\rg$. Hence $\lg a,b\rg\rtimes \lg c,d\rg\unlhd \os$ stabilizes $\lg u\rg$. But on the other hand, applying Clifford theorem we deduce that $\lg a,b\rg\rtimes \lg c,d\rg$ should act irreducibly on $P$, a contradiction to the fact that $\lg a,b\rg\rtimes \lg c,d\rg$ stabilizes $\lg u\rg $.
 \qed

\subsection{Classification of $\MM(G;r,t,\ell)$ with Euler characteristic $-5^4, -7^4$ or $-13^4$}

In this section, we will give a characterization of those regular maps $\MM(G;r,t,\ell)$ with Euler characteristic $-5^4, -7^4$ or $-13^4$ in terms of certain parameters. 
\vskip 2mm
The following lemma give a preliminary characterization of automorphism groups of these maps.
\begin{lem}\label{chara}
Let $\MM(G;r,t,\ell)$ be a regular map with Euler characteristic $-p^4$ for $p\in \{5, 7, 13\}$. Then $O_p(G)\cong\mathbb{Z}^3_p$ is a minimal normal subgroup of $G$, and one of the following holds:
\begin{enumerate}
\item[\rm (1)] $p=5$, $\MM$ has type $\{4,6\}$ and $G/O_5(G)\cong \PGL(2,5)$;
\item[\rm (2)] $p=7$, $\MM$ has type $\{3,8\}$ and $G/O_7(G)\cong \PGL(2,7)$;;
\item[\rm (3)] $p=13$, $\MM$ has type $\{3,7\}$ and $G/O_{13}(G)\cong \PSL(2,13)$;.
\end{enumerate}
\end{lem}
\demo Denote $\MM(G;r,t,\ell)$ by $\MM$. Let $P\in \Syl_p(G)$ and $\{x,y\}$ be the type of $\MM$. According to Proposition \ref{sylow-p}, we deduce that $|P|\geq p^4$. By Theorem \ref{divide} and its proof, we have $p\nmid xy$ and $|P:O_p(G)|\leq p$. Thus the Euler characteristic of the quotient map of $\MM$ induced by $O_p(G)$ is $-\frac{p^4}{|O_p(G)|}$. This together with $|P:O_p(G)|\leq p$ implies that $|O_p(G)|=p^3$ and $|P|=p^4$. By Proposition \ref{cla}, we have three cases:$\{x, y\}=\{4, 6\}$ and $G/O_p(G)\cong \PGL(2,5)$; $\{x, y\}=\{3, 8\}$ and $G/O_p(G)\cong \PGL(2,7)$; or $\{x, y\}=\{3, 7\}$ $G/O_p(G)\cong\PSL(2,13)$. In particular, $G$ is insolvable. 

We show that $O_p(G)$ is a minimal normal subgroup of $G$. If this were proved, then $O_p(G)\cong\mathbb{Z}^3_p$ would follow immediately. Suppose for the contrary there exists a normal subgroup $N$ of $G$ such that $1<N<O_p(G)$. Let $\overline{\MM}$ be the quotient map of $\MM$ induced by $N$. Then $\overline{\MM}$ has type $\{4,6\}$ ( or $\{3,8\}$ or $\{3,7\}$) and $\chi(\overline{\MM})=-\frac{p^4}{|N|}$. Assume $|N|=p^2$. Then if follows from Proposition \ref{Eule} that the type of such a quotient map is $\{3,13\}$, a contradiction. Assume $|N|=p$. Then it follows from Proposition \ref{p^3} that $G/N$ is solvable. This implies that $G$ is solvable, a contradiction to the fact that $G$ is insolvable, as desired. \qed
\begin{rem}
By Lemma \ref{chara}, the automorphism group of a regular map $\MM:=\MM(G;r,t,\ell)$ with $\chi_\MM=-5^4$ (resp. $-7^4$ or $-13^4$) is an extension of $\PGL(2,5)$ (resp. $\PGL(2,7)$ or $\PSL(2,13)$) by one of its $3$-dimensional irreducible module. Note that for a given group $G$ and an $G$-module $M$, the elements in the degree two cohomology group $H^2(G,M)$ are in one-to-one correspondence to the group extensions of $G$ by $M$, where the element $0$ corresponds to a split extension of $G$ by $M$. The group $G=\PGL(2,5)$ (resp. $\PGL(2,7))$ has two $3$-dimension irreducible modules, and exactly one of which, say $M$ (resp. $N$), satisfies that $H^2(G,M)=0$ (resp. $H^2(G, N)=0$). 
\end{rem}

Given a non-orientable regular map $\MM(G;r,t,\ell)$, recall that the {\it reduced presentation} of $G$ is referred to as the presentation of $G$ given by relations in terms of $tr$ and $r\ell$. Let $\MM(\PGL(2,5);r_1,t_1,\ell_1)$, $\MM(\PGL(2,7);r_2,t_2,\ell_2)$ and $\MM(\PSL(2,7);r_3,t_3,\ell_3)$ be regular maps of types $\{4,6\}, \{8,3\}$ and $\{7,3\}$, respectively. Then the reduced presentations of $\PGL(2,5), \PGL(2,7)$ and $\PSL(2,13)$ are given as follows:
$$\PGL(2,5)=\lg R,S\mid R^4=S^6=(RS)^2=(S^{-1}R)^3=1\rg,$$
$$(\rm Pre.1)~~~\PGL(2,7)=\lg R,S\mid R^8=S^3=(RS)^2=(SR^{-2})^4=1\rg,$$
$$(\rm Pre.2)~~~\PGL(2,7)=\lg R,S\mid R^8=S^3=(RS)^2=((SR^{-2})^2R^{-2})^2=1\rg,$$
$$\PSL(2,13)=\lg R,S\mid R^7=S^3=(RS)^2=((SR^{-2})^4SR^3)^2=1\rg.$$
Now, let $\MM(G;r,t,\ell)$ be a regular map Euler characteristic $-p^4$ for $p\in \{5, 7, 13\}$ and let $\overline{\MM}$ be its quotient map induced by $O_p(G)$. Then by Lemma \ref{chara}, we can view $O_p(G)$ as an irreducible $G/O_p(G)$-module and therefore we have seven cases: 

$(1)$ $p=5$ and $H^2(G/O_5(G),O_5(G))= 0$; 

$(2)$ $p=5$ and $H^2(G/O_5(G),O_5(G))\neq 0$; 

$(3)$ $p=7$, $H^2(G/O_7(G),O_7(G))=0$ and the reduced presentation for $\Aut(\overline{\MM})$ is $(\rm Pre.1)$; 

$(4)$ $p=7$, $H^2(G/O_7(G),O_7(G))\neq 0$ and the reduced presentation for $\Aut(\overline{\MM})$ is $(\rm Pre.1)$;

$(5)$ $p=7$, $H^2(G/O_7(G),O_7(G))=0$ and the reduced presentation for $\Aut(\overline{\MM})$ is $(\rm Pre.2)$;

$(6)$ $p=7$, $H^2(G/O_7(G),O_7(G))\neq 0$ and the reduced presentation for $\Aut(\overline{\MM})$ is $(\rm Pre.1)$;

$(7)$ $p=13$. 

According to these seven cases, we define the following seven groups and their related regular maps:
$$\begin{array}{lll}G_1(a_1,b_1,c_1):=&\lg w,s,z|w^4=s^6=z^5=(ws)^2=1, (s^{-1}w)^3=z^{a_1}(z^{s^2})^{b_1}(z^{s^{4}})^{c_1}, \\
& z^w=z^3(z^{s^4})^4,z^{s^2w}=(z^{s^2})^2(z^{s^4})^3,z^{s^4w}=(z^{s^4})^4,z^s=z z^{s^2}(z^{s^{4}})^2,\\
& z^{s^3}=z^2 z^{s^2} z^{s^{4}},z^{s^5}=z(z^{s^2})^2 z^{s^{4}},(z,z^{s^2})=(z,z^{s^4})=1\rg.\end{array}$$where $(a_1,b_1,c_1)\in \mathbb{F}^3_5$ and
$$\begin{array}{lll}G_2(a_2,b_2,c_2):=&\lg w,s,z|w^4=s^6=z^5=(ws)^2=1, (s^{-1}w)^3=z^{a_2}(z^{s^2})^{b_2}(z^{s^{4}})^{c_2}, \\
& z^w=z(z^{s^2})^2 z^{s^4},z^{s^2w}=z^4 z^{s^2}) ,z^{s^4w}=z^{s^2}(z^{s^4})^4,z^s=z^4(z^{s^2})^4(z^{s^{4}})^3,\\
& z^{s^3}=z^3(z^{s^2})^4(z^{s^{4}})^4,z^{s^5}=z^{s^2},(z,z^{s^2})=(z,z^{s^4})=1\rg.\end{array}$$where $(a_2,b_2,c_2)\in \mathbb{F}^3_5$ and$$\begin{array}{lll}G_3(a_3,b_3,c_3):=&\lg w,s,z|w^8=s^3=(ws)^2=z^7=1,(sw^{-2})^4=z^{a_3}(z^s)^{b_3}(z^{s^{-1}})^{c_3}, \\
& z^w=z(z^s)^3(z^{s^{-1}})^6, z^{sw}=z^6z^s(z^{s^{-1}})^4, z^{s^{-1}w}=z^3(z^s)^3 z^{s^{-1}},\\
& (z,z^s)=(z,z^{s^{-1}})=1\rg.\end{array}$$ where $(a_3,b_3,c_3)\in \mathbb{F}^3_7$ and
$$\begin{array}{lll}G_4(a_4,b_4,c_4):=&\lg w,s,z|w^8=s^3=(ws)^2=z^7=1,(sw^{-2})^4=z^{a_4}(z^s)^{b_4}(z^{s^{-1}})^{c_4}, \\
& z^w=z^3(z^s)^3(z^{s^{-1}})^5, z^{sw}=z^4(z^s)^5(z^{s^{-1}})^3, z^{s^{-1}w}=z^5(z^s)^4(z^{s^{-1}})^3,\\
& (z,z^s)=(z,z^{s^{-1}})=1\rg.\end{array}$$ where $(a_4,b_4,c_4)\in \mathbb{F}^3_7$ and
and
$$\begin{array}{lll}G_5(a_5,b_5,c_5):=&\lg w,s,z|w^8=s^3=(ws)^2=z^7=1,((sw^{-2})^2w^{-2})^2=z^{a_5}(z^s)^{b_5}(z^{s^{-1}})^{c_5}, \\
&z^w=z^{s^{-1}}, z^{sw}=z^6(z^s)^2(z^{s^{-1}})^2, z^{s^{-1}w}=z^s,\\
& (z,z^s)=(z,z^{s^{-1}})=1\rg.\end{array}$$ where $(a_5,b_5,c_5)\in \mathbb{F}^3_7$ and
$$\begin{array}{lll}G_6(a_6,b_6,c_6):=&\lg w,s,z|w^8=s^3=(ws)^2=z^7=1,((sw^{-2})^2w^{-2})^2=z^{a_6}(z^s)^{b_6}(z^{s^{-1}})^{c_6}, \\
&z^w=(z^s)^5(z^{s^{-1}})^2, z^{sw}=z^6, z^{s^{-1}w}=(z^s)^5(z^{s^{-1}})^5,\\
& (z,z^s)=(z,z^{s^{-1}})=1\rg.\end{array}$$ where $(a_6,b_6,c_6)\in \mathbb{F}^3_7$ and
$$\begin{array}{lll}G_7(a_7,b_7,c_7):=&\lg w,s,z|w^7=s^3=(ws)^2=z^{13}=1,z^w=z^9(z^s)^5(z^{s^{-1}})^2,z^{sw}=z^{12},\\ 
& z^{s^{-1}w}=z^6(z^s)^{11}(z^{s^{-1}})^2,(z,z^s)=(z,z^{s^{-1}})=1,\\ 
& ((sw^{-2})^4sw^3)^2=z^{a_7}(z^s)^{b_7}(z^{s^{-1}})^{c_7}\rg\end{array}$$where $(a_7,b_7,c_7)\in \mathbb{F}^3_{13}$.
Clearly, for $1\leq i\leq 7$ and any triple $(a_i,b_i,c_i)\neq (0,0,0)$, the group $G_i(a_i,b_i,c_i)$ can be generated by $r$ and $s$ and has trivial center. For any $1\leq i\leq 7$, it follows from that \cite[Theorem 7.7]{NeS} that there is (up to isomorphism and duality) only one regular map of type $\{4, 6\}$ for $i=1,2$ (resp. $\{3, 8\}$ for $i=3, 4, 5,6,$ or $\{3, 7\}$ for $i=7$) with automorphism group $G_i(a_i, b_i, c_i)$, and we denote it by $\MM_i(a_i, b_i, c_i)$.

For convenience, set $p_1=p_2=5, p_3=p_4=p_5=p_6=7$ and $p_7=13$. It is easy to that $\MM_i(a_i,b_i,c_i)\cong \MM_i(ja_i,jb_i,jc_i)$ for $1\leq i\leq 7$ and $1\leq j\leq p_i-1$. Thus we could say that $\MM_i(a_i,b_i,c_i)$ corresponds to a projective point $(a_i,b_i,c_i)$ in  ${\rm PG}(2,p_i)$ for $1\leq i\leq 7$. We say that $(a_i,b_i,c_i)\neq (0,0,0)$ is an {\it admissible point} in ${\rm PG}(2,\mathbb{F}_{p_i})$ if $\MM_i(a_i,b_i,c_i)$ is of Euler characteristic $-p^4_i$, that is, $|G_i(a_i,b_i,c_i)|=\frac{4x_iy_ip^4_i}{x_iy_i-2x_i-2y_i}$ where $\{x_i,y_i\}$ is the type of $\MM_i(a_i,b_i,c_i)$. {\bf Thus,to determine whether a point $(a_i,b_i,c_i)\neq (0,0,0)$ is an admissible point or not, it suffices to check the order of $G_i(a_i,b_i,c_i)$.}
\vskip 2mm
For $1\leq i\leq 7$, let $N_i(a_i,b_i,c_i)$ be the normal closure of $z$ in $G_i(a_i,b_i,c_i)$, and let $\Omega_i$ be the set of regular maps $\MM(G;r,t,\ell)$ (up to isomorphism and duality) with Euler characteristic $-p^4_i$ for which $G$ satisfies the case (i) as above. Let $X_i$ be the matrix in $GL(3, p_i)$ given by the action of $w$ on $N_i(a_i,b_i,c_i)$ with respect to the basis $z, z^s, z^{s^{-1}}$ (or $z, z^{s^2}, z^{s^4}$ when $i=1$), and let $Y_i$ be the permutation matrix in $GL(3, p_i)$ corresponding to the permutation $(123)$. Set $M_i:=Y_iX^{-2}_i$ for $i=1,2$, $M_i:=Y_iX^{-2}_i$ for $i=3,4$, $M_i:=(Y_iX^{-2}_i)^2X^{-2}_i$ for $i=5, 6$, and $M_i:=(Y_iX^{-2}_i)^4Y_iX^3_i$ for $i=7$. Then we have the following result:
\begin{prop}\label{admis}
Let $1\leq i\leq 7$. Then any admissible point $(a_i, b_i,c_i)$ is an eigenvector of eigenvalue $1$ of the matrix $M_i$ and regular maps in $\Omega_i$ are in one-to-one correspondence to the admissible points in ${\rm PG}(2,\mathbb{F}_{p_i})$.
\end{prop}
\demo Since the proofs for the seven cases are completely similar, we just take the case $i=1$ for proof. Let $\MM:=\MM(G;r,t,\ell)\in \Omega_1$ be a regular map with Euler characteristic $-5^4$ and set $R:=rt$ and $S:=r\ell$. Then by Lemma \ref{chara} we deduce that $\MM$ has type $\{4,6\}$, $O_5(G)\cong \mathbb{Z}^3_5$ and $\og:=G/O_5(G)\cong \PGL(2,5)$. This implies that $\MM(\og;\olr, \olt,\oll)$ is a non-orientable regular map of type $\{4,6\}$. By the presentation of $\PGL(2,5)$, we see that $\overline{R}$ and $\overline{S}$ satisfy that $(\overline{R}\overline{S}^{-2})^4=1$, that is, $(RS^{-2})^4\in O_5(G)$. If $(RS^{-2})^4=1$, then it follows from $G=\lg R,S\rg$ and $R^4=S^6=(RS)^2=(RS^{-2})^4=1$ that $5^3\cdot 120=|G|\leq |\PGL(2,5)|=120$, a contradiction. Thus $(RS^{-2})^4\neq 1$. By assumption, the exact sequence $1\ra O_5(G)\ra G\ra G/O_5(G)\ra 1$ is non-split. Thus $H^2(\og, O_5(G))\neq 0$, where $O_5(G)$ is viewed as an $\og$-module induced by the conjugation of $G$ on $O_5(G)$. Consider the 3-dimensional representation, say $\rho:\og\ra \GL(3,5)$, afforded by $O_5(G)$. Then $\rho(S^2)$ is similar to the permutation matrix given by the permutation $(123)$, with transformation matrix $X$. From this, we see that these are only three elements $z, z^{S^2}, z^{S^4}\in O_5(G)$ (up to power) such that $O_5(G)=\lg z, z^{S^2}, z^{S^4}\rg $, and the conjugation of $R$ on $O_5(G)$ is given by the matrix $X^{-1}\rho(\overline{R})X$. The matrix $X^{-1}\rho(\overline{R})X$ is fixed if and only if $z$ is fixed up to power. Now, fix $z$ such that $X^{-1}\rho(\overline{R})X$ corresponding to the action of $w$ on $z, z^{s^2}$ and $z^{s^4}$ in the representation of $G_1(a_1,b_1,c_1)$. Since $(RS^{-2})^4\in O_5(G)$, there exist non-negative integers $a, b$ and $c$ such that $(RS^{-2})^4=z^a(z^{S^2})^b(z^{S^{4}})^c$. Clearly, $(a,b,c)\neq (0,0,0)$ as $(RS^{-2})^4\neq 1$. Hence the presentation of $G$ in terms of $R,S$ and $z$ is the same as that of $G_1(a,b,c)$ in terms of $r, s$ and $z$. This implies that $\MM$ is isomorphic to $\MM_1(a,b,c)$ for some admissible point $(a,b,c)$. Since by definition $M_1$ is the matrix given by the action of $RS^{-2}$ on $O_5(G)$ with respect to the basis $z, z^{s^2}$ and $z^{s^4}$ and $RS^{-2}$ centralizes its any power, it follows that $(a,b,c)$ is an eigenvector of eigenvalue $1$ of the matrix $M_1$.

To complete the proof, it suffices to show that $\MM_1(a_1,b_1,c_1)$ and $\MM_1(a'_1,b'_1,c'_1)$ are non-isomorphic for any two distinct admissible points $(a_1, b_1, c_1)$ and $(a'_1,b'_1,c'_1)$ in ${\rm PG}(2, \mathbb{F}_5)$. Suppose for the contrary, $\MM_1(a_1,b_1,c_1)$ and $\MM_1(a'_1,b'_1,c'_1)$ are isomorphic for two distinct admissible points $(a_1, b_1, c_1)$ and $(a'_1,b'_1,c'_1)$. Let $G_1:=G_1(a_1,b_1,c_1)$ and $G'_1:=G_1(a'_1,b'_1,c'_1)$, and let $\MM(G_1;r_1,t_1,\ell_1)$ and $\MM(G'_1;r'_1,t'_1,\ell'_1)$ denote the regular maps $\MM_1(a_1,b_1,c_1)$ and $\MM_1(a'_1,b'_1,c'_1)$, respectively. Then by assumption there is an isomorphism $\varphi: G_1\ra G'_1$ such that $\varphi(r)=r', \varphi(t)=t'$ and $\varphi(\ell)=\ell'$. Set $R_1:=r_1t_1, S_1:=r_1\ell_1, R'_1:=r'_1t'_1$ and $S'_1:=r'_1\ell'_1$, and denote by $z_1$ and $z'_1$ the unique (up to power) element (given by rule in previous paragraph) in $O_5(G_1)$ and $O_5(G'_1)$, respectively. Then $R_1, S_1$ and $z_1$ (resp. $R'_1, S'_1$ and $z'_1$) satisfy the same relations as for $r, s$ and $z$ in the presentation of $G_1(a_1,b_1,c_1)$ (resp. $G_1(a'_1,b'_1,c'_1)$). Since $z_1$ (resp. $z'_1$) is uniquely determined (up to power) by $S_1$ and $R_1$ (resp. $S'_1$ and $R'_1$) by the rule that $O_5(G_1)=\lg z_1, {z_1}^{S^2_1}, {z_1}^{S^4_1}\rg$ (resp. $O_5(G'_1)=\lg z'_1, {z'_1}^{{S'_1}^2}, {z'_1}^{{S'_1}^4}\rg $) and the conjugation of $R_1$ (resp. $R'_1$) on $O_5(G_1)$ (resp. $O_5(G'_1)$) is given by a fixed matrix with respect to this basis, it follows that $\varphi(z_1)=(z'_1)^m$ for some integer $1\leq m\leq 4$. From $\varphi((S_1R^{-2}_1)^3)= (\varphi(S_1)\varphi(R_1)^{-2})^3=(S'_1{R'_1}^{-2})^3$, we see that $m(a_1,b_1,c_1)=(a'_1,b'_1,c'_1)$. This is a contradiction to the assumption that $(a_1, b_1, c_1)$ and $(a'_1,b'_1,c'_1)$ represent distinct points in ${\rm PG}(2,\mathbb{F}_5)$. \qed
\vskip 2mm
By Proposition \ref{admis} and direct computation with MAGMA, we arrive at the following goal:
\begin{theorem}\label{cha-exc}
A regular map $\MM=\MM(G;r,t,\ell)$ is of Euler characteristic $-p^4$ for $p\in \{5, 7, 13\}$ if and only if one of the following holds:
\begin{enumerate}
\item[\rm (1)] $p=5$, $\MM$ is isomorphic to $\MM_1(1,3,4)$ or its dual;
\item[\rm (2)] $p=7$, $\MM$ is isomorphic to $\MM_4(1,1,4)$ or $\MM_6(1,1,3)$ or their dual;
\item[\rm (3)] $p=13$, $\MM$ is isomorphic to $\MM_7(1,7,6)$ or its dual.
\end{enumerate}
\end{theorem}
\vskip 2mm
\noindent {\bf Acknowledgement:} The first author is partially supported by the National Natural Science Foundation of China (Grant Nos. 12350710787, 12571016), the second author is partially supported by the National Natural Science Foundation of China (Grant No. 12471332).

\vskip 2mm
\noindent {\bf Data Availability Statement}
Data sharing is not applicable to this article as no new data were created or analyzed in this study.
\vskip 2mm
\noindent {\bf Disclosure Statement}
All authors do not have any relevant financial or non-financial competing interests.

\end{document}